\documentclass[12pt,reqno]{amsart}
\usepackage[letterpaper,margin=1in]{geometry}
\usepackage[T1]{fontenc}
\usepackage{placeins}
\usepackage{capt-of}
\usepackage{xcolor}
\usepackage{lmodern,microtype,amsmath,amssymb,mathtools,mathrsfs}
\microtypesetup{expansion=false}
\usepackage{amsthm}
\usepackage[pagebackref]{hyperref}
\hypersetup{hidelinks}
\usepackage{booktabs,longtable,array,enumitem}
\usepackage{tikz}
\usetikzlibrary{arrows.meta,calc,positioning,patterns,decorations.pathreplacing}
\numberwithin{equation}{section}
\setlist[enumerate]{leftmargin=*,itemsep=3pt,topsep=4pt}
\newtheorem{theorem}{Theorem}[section]
\newtheorem{proposition}[theorem]{Proposition}
\newtheorem{lemma}[theorem]{Lemma}
\newtheorem{corollary}[theorem]{Corollary}
\theoremstyle{definition}
\newtheorem{definition}[theorem]{Definition}
\newtheorem*{example}{Example}
\theoremstyle{remark}
\newtheorem*{remark}{Remark}
\theoremstyle{plain}
\newtheorem{conjecture}{Conjecture}
\newcommand{\NN}{\mathbb Z_{\ge0}}
\newcommand{\ZZ}{\mathbb Z}
\newcommand{\QQ}{\mathbb Q}
\newcommand{\one}{\mathbf1}
\newcommand{\nn}{\mathbf n}
\newcommand{\cc}{\mathbf c}

\newcommand{\qbin}[2]{\genfrac{[}{]}{0pt}{}{#1}{#2}_{q}}
\newcommand{\mult}[2]{\left[\begin{matrix}#1\\#2\end{matrix}\right]_{q;G}}
\DeclareMathOperator{\area}{area}
\DeclareMathOperator{\arm}{arm}
\DeclareMathOperator{\leg}{leg}
\DeclareMathOperator{\dinv}{dinv}
\DeclareMathOperator{\tdinv}{tdinv}
\DeclareMathOperator{\dist}{dist}

\DeclareMathOperator{\ides}{ides}

\newcommand{\Qop}{\mathsf Q}
\newcommand{\Dop}{\mathsf D}
\newcommand{\Cop}{\mathsf C}
\DeclarePairedDelimiter\abs{\lvert}{\rvert}

\newcommand{\Hseries}{\mathscr H}
\newcommand{\Aseries}{\mathscr A}
\newcommand{\Dseries}{\mathscr D}
\newcommand{\Bseries}{\mathscr B}
\newcommand{\Pseries}{\mathscr P}
\newcommand{\Rpaths}{\mathcal R_{a,b}^{(N)}}
\newcommand{\Ggraph}{\mathcal G_{a,b}}

\title{Rogers--Ramanujan identities from the geometry of $X^a=Y^b$}

\author{Yifeng Huang}
\address{Department of Mathematics, Emory University, 400 Dowman Drive, Atlanta, GA 30322}
\email{yifeng.huang@emory.edu}

\author{Kenny Lau}
\address{Axiom Math, 124 University Avenue, Palo Alto, CA 94301}
\email{kenny@axiommath.ai}
\author{Ken Ono}
\address{Axiom Math, 124 University Avenue, Palo Alto, CA 94301}
\email{ken@axiommath.ai}

\date{}
\begin{document}
\begin{abstract}
We prove the conjecture of Huang, Jiang, and Oblomkov (HJO)
giving a geometric extension of the Rogers--Ramanujan and
Andrews--Gordon identities for every torus-knot singularity
\(X^a=Y^b\) with coprime \(1<a<b\).
For a prime power \(q\), let
\(\mathcal{NC}_n^{a,b}(\mathbb F_q)\) denote the set of pairs
of commuting nilpotent \(n\times n\) matrices \((A,B)\) over
\(\mathbb F_q\) satisfying \(A^a=B^b\).
We establish the threefold equality between their normalized
counts, the HJO \(q\)-series \(Z_{a,b}\), and the explicit
infinite product \(P_{a,b}\):
\[
\underbrace{\vphantom{\Bigg|}
\prod_{m\geq1}(1-q^{-m})
\Biggl(\sum_{n=0}^{\infty}
\frac{\lvert\mathcal{NC}_n^{a,b}(\mathbb F_q)\rvert}
{\lvert\operatorname{GL}_n(\mathbb F_q)\rvert}\Biggr)
}_{\text{point count}}
=
\underbrace{\vphantom{\Bigg|}Z_{a,b}(q^{-1})
}_{\text{\(q\)-series}}
=
\underbrace{\vphantom{\Bigg|}P_{a,b}(q^{-1})
}_{\text{infinite product}}.
\]
Our main result is a stronger finite identity: the rank \(N\)
HJO sum equals \((q;q)_N\) times the generating function for
balanced cylindric partitions with entries bounded by \(N\).
Taking \(N\to\infty\) yields the HJO conjecture.
The proof combines the compositional rational shuffle theorem
of Bergeron--Garsia--Leven--Xin and Mellit with a multiplicativity
theorem for slope operators and a determinantal model for
bounded cylindric partitions, linked by a common
\(q\)-difference equation.
The finite identity and the HJO conjecture have been formalized unconditionally
in Lean by AxiomProver.
\end{abstract}

\maketitle

\section{Introduction}

We begin with the Rogers--Ramanujan identities, which are the classical identities that our main result extends.
For an indeterminate $q$, write
\[
 (x;q)_n:=\prod_{j=0}^{n-1}(1-xq^j),\qquad
 (x;q)_\infty:=\prod_{j\ge0}(1-xq^j),\qquad
 (x_1,\ldots,x_r;q)_n:=\prod_{i=1}^r(x_i;q)_n.
\]
The identities are
\begin{align}
 \sum_{n\ge0}\frac{q^{n^2}}{(q;q)_n}
 &=\frac1{(q,q^4;q^5)_\infty},\label{eq:RRfirst}\\
 \sum_{n\ge0}\frac{q^{n^2+n}}{(q;q)_n}
 &=\frac1{(q^2,q^3;q^5)_\infty}.\label{eq:RRsecond}
\end{align}
An empty product is one. Infinite products with positive powers of $q$ are interpreted coefficientwise. These identities may also be read analytically for $|q|<1$. Their sum sides count partitions whose successive parts differ by at least two, with smallest part at least one or at least two, respectively. Their product sides count partitions with parts congruent to $1,4$ or $2,3$ modulo five, respectively. This equality of two apparently different combinatorial counts is the basic phenomenon that we generalize (for example, see Rogers~\cite{Rogers} and Andrews~\cite[Chapter~7]{AndrewsBook}).

Gordon's partition theorem~\cite{Gordon} and Andrews's analytic identities~\cite[pp.~4082--4083]{Andrews} extend the modulus five in \eqref{eq:RRfirst}--\eqref{eq:RRsecond} to every odd modulus. Namely, for integers $k\ge2$ and $1\le i\le k$, we have
\begin{equation}\label{eq:AG}
 \sum_{n_1\ge\cdots\ge n_{k-1}\ge0}
 \frac{q^{n_1^2+\cdots+n_{k-1}^2+n_i+\cdots+n_{k-1}}}
 {(q;q)_{n_1-n_2}\cdots(q;q)_{n_{k-2}-n_{k-1}}(q;q)_{n_{k-1}}}
 =\prod_{\substack{m\ge1\\m\not\equiv0,\,\pm i\pmod{2k+1}}}
 (1-q^m)^{-1}.
\end{equation}
The linear sum in the exponent is empty when $i=k$. In particular, $k=2$, $i=2$ and $i=1$ recover the Rogers--Ramanujan identities.
Unlike an unrestricted partition generating function, the multiple sum in \eqref{eq:AG} includes a quadratic weight and a sequence of inequalities. These features will reappear in the sums considered below.

\subsection{The HJO conjecture}
For relatively prime integers $1<a<b$, we
let
\begin{equation}\label{eq:parameters}
 d:=a+b,\qquad \delta:=\frac{(a-1)(b-1)}2,\qquad f:=ab-a-b.
\end{equation}
Write $\NN:=\{0,1,2,\ldots\}$. The {\it numerical semigroup}
\(
 \Gamma:=\langle a,b\rangle=\{ua+vb:u,v\in\NN\}
\)
has gap set $G:=\NN\setminus\Gamma$. It has $\delta$ gaps and largest gap $f$. We make $G$ a poset by
\[
 g\preceq h\quad\Longleftrightarrow\quad h-g\in\Gamma.
\]
We note that $f$ is the maximum for this partial order. The Sylvester coordinates used in Section~\ref{sec:definitions} make these facts explicit.

A {\it gap vector} is a tuple $\nn=(n_g)_{g\in G}$. We let
\[
 \mathcal C:=\{\nn\in\NN^G:n_g\le n_h\text{ whenever }g\preceq h\}.
\]
A gap vector is extended to a vector on $\ZZ$ by setting $n_j=0$ for $j\not \in G$. Since $f$ is the maximum gap, $n_f$ is the largest coordinate of
a vector in $\mathcal C$.  Write $(q)_m:=(q;q)_m$.  For a gap vector $\nn$ and an integer $N\ge0$, we let

\begin{align}
 K(s)&:=\one_{s\ge0}-\one_{s\ge a}-\one_{s\ge b}+\one_{s\ge d},\label{eq:K}\\
 Q(\nn)&:=\sum_{g,h\in G}K(h-g)n_gn_h,\label{eq:Q}\\
 P_G(\nn)&:=\prod_{g\in G}
 \frac{(q)_{n_g-n_{g-a-b}}}{(q)_{n_g-n_{g-a}}(q)_{n_g-n_{g-b}}},\label{eq:PG}\\
 \mult N\nn&:=\frac{(q)_N}{(q)_{N-n_f}}P_G(\nn)
 \qquad(\nn\in\mathcal C,\ n_f\le N).\label{eq:mult}
\end{align}
We call \eqref{eq:mult} a \emph{generalized Gaussian multinomial}.
Its factorization into ordinary Gaussian polynomials in
Section~\ref{sec:definitions} shows that it belongs to $\NN[q]$.

For an indeterminate $t$ and an integer $N\ge0$, the first author, Jiang, and Oblomkov (HJO)~
\cite[(2.1)--(2.5)]{HJO} define
\begin{align}
 &N_{a,b;N}(q,t)
:=\sum_{\substack{\nn\in\mathcal C\\n_f\le N}}
       \mult N\nn q^{Q(\nn)}t^{|\nn|},
       \qquad |\nn|:=\sum_{g\in G}n_g,\label{eq:HN}\\
 F_N(q)&:=N_{a,b;N}(q,1),\qquad
 Z_{a,b}(q):=N_{a,b;\infty}(q,1):=\sum_{\nn\in\mathcal C}P_G(\nn)q^{Q(\nn)}.
 \label{eq:Z}
\end{align}
The first author proved the bound~\cite[Theorem~1.3]{Huang}
\begin{equation}\label{eq:coercivity}
 Q(\nn)\ge n_f^2/\delta\qquad(\nn\in\mathcal C).
\end{equation}
This is an effective positive definiteness theorem on the cone $\mathcal C$, not on the ambient vector space.
It shows that the infinite sum in \eqref{eq:Z} is well defined coefficientwise and that $F_N$ tends to $Z_{a,b}$ in that sense.

The geometric source of these expressions is the singularity
$R=\Bbbk[[T^a,T^b]]$, with normalization $\widetilde R=\Bbbk[[T]]$.
The punctual Quot scheme $\operatorname{Quot}_m^R(\widetilde R^N)$
parametrizes $R$-submodules of $\widetilde R^N$ of colength $m$.
HJO gives an affine-cell decomposition and expresses its motivic
generating function through \eqref{eq:HN}.  Therefore, $N$ is the module rank in the geometric interpretation.
For the arithmetic interpretation, let $q$ be a prime power and set
\[
 \mathcal{NC}_m^{a,b}(\mathbb F_q)
 :=\{(A,B)\in\operatorname{Mat}_m(\mathbb F_q)^2:
       A,B\text{ nilpotent},\ AB=BA,\ A^a=B^b\}.
\]
Then \cite[Theorem~1.4]{HJO} gives
\begin{equation}\label{eq:matrixvolume}
 \sum_{m\ge0}
 \frac{|\mathcal{NC}_m^{a,b}(\mathbb F_q)|}{|\mathrm{GL}_m(\mathbb F_q)|}
 =\frac{Z_{a,b}(q^{-1})}{(q^{-1};q^{-1})_\infty}.
\end{equation}
The left side is the groupoid volume of the category of finite modules over
$\mathbb F_q[[T^a,T^b]]$.
Only in this finite-field interpretation does $q$ denote a prime power.
In the formal identities below, it is again an indeterminate.
These interpretations are established in HJO independently of the
product evaluation.

HJO conjectured an infinite product evaluation of $Z_{a,b}(q)$ with nonnegative exponents periodic modulo $d=a+b$.  For an integer $x$, put $\dist(x,d\ZZ):=\min_{j\in\ZZ}|x-dj|$. Their prediction is the following.

\begin{conjecture}[{HJO \cite[Conjecture~1.7]{HJO}}]\label{conj:HJO}
For relatively prime integers $1<a<b$, we have
\begin{equation}\label{eq:HJOconj}
 Z_{a,b}(q)=P_{a,b}(q):=\prod_{m\ge1}(1-q^m)^{-\rho_{a,b}(m)},
\end{equation}
where 
$$
 \rho_{a,b}(m):=\min\{a,b,\dist(am,d\ZZ)\}-1+\one_{d\mid m}.
$$
\end{conjecture}

\begin{remark} The nonnegative exponents in the conjecture are periodic with period $d$. We note that they also vanish on the zero residue class. Therefore, the conjecture predicts another equality between a constrained quadratic sum and a product with congruence conditions. Moreover, these infinite products are modular functions up to a rational power of $q$.
\end{remark}

\begin{example}[Rogers--Ramanujan and Andrews--Gordon]\label{ex:classical}
Let $a=2$ and $b=2m+1$ with $m\ge1$, so that $d=2m+3$. The gap set is the chain $G=\{1,3,\ldots,2m-1\}$, in which each gap is obtained from its predecessor by adding $a=2$. Off-$G$ coordinates vanish, so \eqref{eq:K}--\eqref{eq:PG} specialize to
\[
 Q(\nn)=\sum_{g\in G}n_g^2,\qquad
 P_G(\nn)=\prod_{g\in G}\frac{1}{(q)_{n_g-n_{g-2}}},
\]
and the cone $\mathcal C$ is the chain of inequalities $n_1\le n_3\le\cdots\le n_{2m-1}$. The change of variables $N_j:=n_{2(m-j)+1}$, for $1\le j\le m$, therefore identifies $Z_{2,2m+1}(q)$ with the left side of \eqref{eq:AG} for $k=i=m+1$, where the linear part of the exponent is empty. On the product side, $\rho_{2,2m+1}$ only takes the values $0$ and $1$, and $\rho_{2,2m+1}(n)=0$ exactly when $2n\equiv0$ or $\pm1\pmod{2m+3}$ (i.e. when $n\equiv0,\pm(m+1)\pmod{2m+3}$). Therefore, \eqref{eq:HJOconj} for $(2,2m+1)$ is precisely the Andrews--Gordon identity \eqref{eq:AG} with $k=i=m+1$. For  $(a,b)=(2,3),$ it is the first Rogers--Ramanujan identity \eqref{eq:RRfirst}. 
\end{example}

Apart from the classical $a=2$ cases of the HJO conjecture, little was known. The authors and Paule \cite{HLOP} proved the $(a,b)=(3, 4), (3, 5), (3,7),$ and $(3,8)$ cases of the conjecture using tricky symbolic $q$-series algebra, and two of the authors \cite{LO} later proved the complete $a=3$ layer of the conjecture. Nothing was known for $a>3.$

\begin{example}[The case $(a,b)=(3,4)$]\label{ex:34}
Here $d=7$ and $G=\{1,2,5\}$, with $1\preceq5$ and $2\preceq5$, while $1$ and $2$ are incomparable. The quadratic form \eqref{eq:Q} is $Q(\nn)=n_1^2+n_2^2+n_5^2+n_1n_2-n_1n_5$, and the conjecture gives the identity
\begin{equation}\label{eq:example34}
 \sum_{\substack{n_1,n_2,n_5\ge0\\ n_1,\,n_2\le n_5}}
 \frac{(q)_{n_5}\,q^{\,n_1^2+n_2^2+n_5^2+n_1n_2-n_1n_5}}
      {(q)_{n_1}(q)_{n_2}(q)_{n_5-n_1}(q)_{n_5-n_2}}
 =\frac{1}{(q,q^{6};q^{7})_\infty^{2}\,(q^{3},q^{4};q^{7})_\infty}.
\end{equation}
This is the smallest case outside the classical family $a=2$. It belongs to the $a=3$ family proved by the second and the third authors (LO)~\cite{LO}. Note the new phenomena visible already here, where the gap poset is no longer a chain, and $Q$ is no longer a sum of squares.
\end{example}

\begin{example}[The case $(a,b)=(4,5)$]\label{ex:45}
Here $d=9$ and $G=\{1,2,3,6,7,11\}$, with cover relations $1,2\preceq6$, $\ 2,3\preceq7$, and $6,7\preceq11$. The quadratic form \eqref{eq:Q} is
\[
\begin{aligned}
 Q(\nn)={}&n_1^2+n_2^2+n_3^2+n_6^2+n_7^2+n_{11}^2
+n_1n_2+n_1n_3+n_2n_3+n_3n_6+n_6n_7\\
 &\ \ \ \ \ -n_1n_6-n_1n_7-n_2n_7-n_3n_{11}-n_6n_{11},
\end{aligned}
\]
and the conjecture gives
\begin{displaymath}
\begin{split}
 \sum_{\nn\in\mathcal C}
 \frac{(q)_{n_6}(q)_{n_7}(q)_{n_{11}-n_2}\,q^{\,Q(\nn)}}
      {(q)_{n_1}(q)_{n_2}(q)_{n_3}
       (q)_{n_6-n_1}(q)_{n_6-n_2}(q)_{n_7-n_2}(q)_{n_7-n_3}
       (q)_{n_{11}-n_6}(q)_{n_{11}-n_7}}\hspace{3.5cm}&\\
 =\frac{1}{(q,q^{8};q^{9})_\infty^{3}\,(q^{3},q^{6};q^{9})_\infty^{2}\,(q^{4},q^{5};q^{9})_\infty}&,
\end{split}
\end{displaymath}
where the sum is over $(n_1,n_2,n_3,n_6,n_7,n_{11})\in\NN^6$ satisfying $n_1,n_2\le n_6$, $\ n_2,n_3\le n_7$ and $n_6,n_7\le n_{11}$. This is the smallest case not covered by the earlier results for $a\le3$, and hence the first identity established by the present paper. We revisit this example in Section~\ref{sec:examples}.
\end{example}

For $a=3$, LO~\cite{LO} first fix a boundary chain of gap
variables and identify the remaining finite sum with a completely
symmetric $q$-supernomial.  Schilling--Warnaar's supernomial theory
\cite{SW,Schilling} transforms this expression into the summands of
Warnaar's $A_2$ Andrews--Gordon identities
\cite[(1.10)--(1.11)]{Warnaar}.  The residues $b\equiv\pm1\pmod3$
correspond to two orders of the tableau blocks.  Restoring the boundary
factors and summing the remaining variables gives the required product. The key point is that LO's proof is a transformation of finite sums, while Warnaar \cite[(1.10)--(1.11)]{Warnaar} supplies the infinite product.

That approach does not easily extend to general $a$. Instead, we connect the HJO finite sums $F_N(q)$ to generating functions
for cylindric partitions with largest entries bounded by $N$. The HJO conjecture is obtained by letting $N\to \infty$. For these finite sums, results of Borodin \cite{Borodin} and Foda--Welsh \cite{FW} on unbounded cylindric partitions supply the infinite product. 

An ordinary partition $\lambda=(\lambda_i)_{i\in \ZZ_{\ge 1}}$ is a weakly decreasing sequence of nonnegative integers with finite sum. Trailing zeros are understood. Its size $|\lambda|$ is the sum of its parts. A {\it profile} is a tuple $\cc:=(c_i)_{i \in \ZZ/a\ZZ}=(c_0,\ldots,c_{a-1})$ of nonnegative integers. The indices are understood cyclically modulo $a$. A {\it cylindric partition} with \emph{outgoing profile} $\cc$ is an $a$-tuple of partitions $\boldsymbol\lambda=(\lambda^{(i)})_{i\in \ZZ/a\ZZ}$ satisfying
\begin{equation}\label{eq:cylinder}
 \lambda^{(i)}_j\ge\lambda^{(i+1)}_{j+c_i}
 \qquad(0\leq i<a,\ j\ge1).
\end{equation}
The partitions $\lambda^{(i)}$ are called the rows of $\boldsymbol\lambda$. The last inequality links the last row back to the first. The volume is $\abs{\boldsymbol{\lambda}}\coloneqq \sum_{i=0}^{a-1}|\lambda^{(i)}|$. Denote by $C_{\cc,\le N}(q)$ its volume generating function when every entry is at most $N$, and by $C_{\cc}(q)$ the same function without this bound. 
Cylindric partitions were introduced by Gessel and Krattenthaler~\cite{GK}.  We use the conventions and product theory of Borodin~\cite[Section~5]{Borodin} and Foda--Welsh~\cite[Sections~2--3]{FW}.

We use the balanced profile
\begin{equation}\label{eq:balanced}
 c_i:=\left\lfloor\frac{(i+1)b}{a}\right\rfloor-
         \left\lfloor\frac{ib}{a}\right\rfloor,
 \qquad0\le i<a.
\end{equation}
Its entries sum to $b$ and each entry is either $\lfloor{b/a}\rfloor$ or $\lceil{b/a}\rceil$. For example, $(a,b)=(5,7)$ gives $\cc=(1,1,2,1,2)$.
Our main result gives a closed expression for the bounded cylindric partition generating function.

\begin{theorem}\label{thm:mainfinite}
For coprime integers $1<a<b$, the profile \eqref{eq:balanced}, and every integer $N\ge0$, we have
\begin{equation}\label{eq:mainfinite}
 N_{a,b;N}(q,1)=(q;q)_N C_{\cc,\le N}(q).
\end{equation}
\end{theorem}

Although the cylindric partition has bounded entries, its rows can have arbitrarily many parts. In particular, $C_{\cc,\le N}(q)$ is generally an infinite series. 
Letting $N\to\infty$ in Theorem~\ref{thm:mainfinite} and invoking the established product evaluation of the unbounded cylindric series settles the conjecture in full.
\begin{corollary}\label{cor:HJO}
 Conjecture~\ref{conj:HJO} is true. 
\end{corollary}
HJO identifies $P_{a,b}(q)$ with the normalized character of the balanced
sector of $W_a(a,a+b)$, labelled by the positive tuple $\cc+\one$
\cite[Proposition~9.3]{HJO}. It then follows from Foda and Welsh \cite[(5), (23), (25), and Sections 3.1–3.3]{FW} that $P_{a,b}(q)=(q;q)_\infty C_\cc(q)$. Hence, Theorem~\ref{thm:mainfinite} implies Corollary~\ref{cor:HJO}.

The significance of Theorem~\ref{thm:mainfinite} goes beyond the HJO conjecture. Our proof identifies $F_N(q)$ with several combinatorial generating functions involving cylindric partitions, prescribed-return composite Dyck paths, and vertex-disjoint cycles of the rank graph (see Corollary~\ref{cor:threemodels}). In particular, these combinatorial generating functions are the same, and \eqref{eq:HN}--\eqref{eq:Z} gives a common analytic formula for them. 

Especially interesting is the consequence of \eqref{eq:mainfinite} that $(q;q)_NC_{\cc,\le N}(q)$ is a polynomial with nonnegative coefficients. For a general $r$-row profile $\cc=(c_0,\dots,c_{r-1})\in \NN^r$ with level $c=c_0+\dots+c_{r-1}>0$, Kur\c{s}ung\"oz--\"Omr\"uuzun Seyrek showed that $(q^r;q^r)_NC_{\cc,\le N}(q)$ is a polynomial with nonnegative coefficients
\cite[(2), arXiv version~2]{KOS}. For a coprime profile, meaning $\gcd(r,c)=1$, Corteel--Warnaar conjectured a stronger statement that $(q;q)_NC_{\cc,\le N}(q)$ is a polynomial with nonnegative coefficients. In fact, what they conjectured is a finer positivity \cite[Conjecture~8.4]{Warnaar}. For a coprime profile $\cc$, define $Q_{m,\cc}(q)$ by
\begin{equation}\label{eq:WarnaarExpansion}
 (q;q)_N C_{\cc,\le N}(q)
   =\sum_{m=0}^N\qbin Nm Q_{m,\cc}(q),
\end{equation}
where $\qbin Nm=(q)_N/((q)_m(q)_{N-m})$ is the ordinary Gaussian
binomial. This is equivalent to the definition \cite[(2.4)]{Warnaar} in the coprime case. By a triangularity argument, $Q_{m,\cc}(q)$ is uniquely defined by this expression, and its polynomiality for all $m$ is equivalent to that of $(q;q)_N C_{\cc,\le N}(q)$ for all $N$. \cite[Conjecture~8.4]{Warnaar} conjectures that $Q_{m,\cc}(q)\in \NN[q]$ and gives a closed form expression for $Q_{m,\cc}(1)$. According to \cite[following Conjecture~8.4]{Warnaar}, the polynomiality and the $q=1$ value is established by unpublished work of Welsh. Therefore, the positivity of $Q_{m,\cc}(q)$ remains the deepest part of the conjecture, and it implies the positivity of $(q;q)_NC_{\cc,\le N}(q)$ by \eqref{eq:WarnaarExpansion}.

Our expression for $N_{a,b;N}(q,1)$ implies that the finer positivity holds for coprime balanced profiles \eqref{eq:balanced}.

\begin{corollary}\label{cor:Warnaar}
For every $m\ge0$ and the balanced profile \eqref{eq:balanced}, we have
\begin{equation}\label{eq:positiveWarnaar}
 Q_{m,\cc}(q)=
 \sum_{\substack{\nn\in\mathcal C\\n_f=m}}
       q^{Q(\nn)}\mult m\nn\ \in\NN[q].
\end{equation}
Thus the positivity assertion of Conjecture~8.4 of \cite{Warnaar} holds for these balanced coprime
profiles.
\end{corollary}
\begin{proof}
    [Proof assuming Theorem~\ref{thm:mainfinite}]
    Compare \eqref{eq:HN}--\eqref{eq:Z} and \eqref{eq:WarnaarExpansion}.
\end{proof}
\begin{remark}
    The $Q_{m,\cc}(1)$ assertion of \cite[Conjecture~8.4]{Warnaar} follows from \eqref{eq:positiveWarnaar} and the word-of-filters model in Section~\ref{sec:return}. We leave it as an exercise for the interested readers.
\end{remark}

\subsection{\texorpdfstring{The first finite Rogers--Ramanujan identity}{The first finite Rogers--Ramanujan identity}}\label{subsec:runninggap}

For $(a,b)=(2,3)$, the gap set is $G=\{1\}$ and $\delta=1$.
A gap vector is a single integer $m$ with $0\le m\le N$. Therefore, \eqref{eq:Q} and \eqref{eq:mult} give
\[
 Q(m)=m^2,\qquad
 \mult Nm=\frac{(q)_N}{(q)_m(q)_{N-m}},\qquad
 N_{2,3;N}(q,t)=\sum_{m=0}^{N}\frac{(q)_N q^{m^2}t^m}{(q)_m(q)_{N-m}}.
\]
This specialization is a finite Gaussian version of the first Rogers--Ramanujan sum. At $N=2$ the three values $m=0,1,2$ give
\begin{equation}\label{eq:exampleF2}
 N_{2,3;2}(q,t)=1+(q+q^2)t+q^4t^2,
 \qquad F_2(q)=1+q+q^2+q^4.
\end{equation}
The balanced profile is $(1,2)$, so Theorem~\ref{thm:mainfinite} asserts in this case that
\[
 C_{(1,2),\le2}(q)=\frac{1+q+q^2+q^4}{(1-q)(1-q^2)}.
\]
Hence, the finite form has the Gaussian numerator while bounding the largest cylinder entry.

\subsection{Organization of the proof}\label{subsec:roadmap}
We prove the finite identity by showing that two generating series satisfy
the same $q$-difference equation.  Its coefficients involve the areas of
rational Dyck paths, which we now define.

Let $\mathcal D_{aN,bN}$ be the set of paths from $(0,0)$ to $(aN,bN)$
with unit north and east steps, weakly above the diagonal.
If $y_x$ is the height of the east step from column $x$ to column $x+1$, put
\begin{equation}\label{eq:introarea}
 \area(P):=\sum_{x=0}^{aN-1}
       \left(y_x-\left\lceil\frac{b(x+1)}a\right\rceil\right),
 \qquad A_N(q):=\sum_{P\in\mathcal D_{aN,bN}}q^{\area(P)}.
\end{equation}
This counts full lattice squares between the path and the diagonal.
Set $A_0(q)=1$ and
\begin{equation}\label{eq:shifts}
 \kappa_N:=(d-1)\binom N2,\qquad \gamma_N:=d\binom N2.
\end{equation}
We then define
\begin{equation}\label{eq:Hseries}
 \begin{aligned}
  \Hseries(z;q)&:=\sum_{N\ge0}q^{\gamma_N}\frac{F_N(q)}{(q)_N}z^N,
  &\Aseries(z;q)&:=\sum_{N\ge0}A_N(q)z^N,\\
  \Bseries(z;q)&:=\sum_{N\ge0}q^{\gamma_N}C_{\cc,\le N}(q)z^N.
 \end{aligned}
\end{equation}
Therefore, Theorem~\ref{thm:mainfinite} is equivalent to $\Hseries=\Bseries$.
We work in $\QQ((q))[[z]]$, the ring of formal power series in $z$
with coefficients that are Laurent series in $q$ with finitely many
negative powers.  Consider the equation
\begin{equation}\label{eq:commonequation}
 \mathscr G(qz;q)=\Aseries(-z;q)\mathscr G(z;q),\qquad
 \mathscr G(0;q)=1.
\end{equation}
Its unique solution is $\mathscr G(z;q)=\sum_{N\ge0}G_N(q)z^N$, where
\begin{equation}\label{eq:solutionrecurrence}
 G_0(q)=1,\qquad
 G_N(q)=\frac{1}{1-q^N}\sum_{k=1}^N(-1)^{k+1}A_k(q)G_{N-k}(q)
 \quad(N\ge1).
\end{equation}
Indeed, coefficient comparison gives this recurrence.  Since $1-q^N$
is invertible, it both determines the coefficients and constructs a
solution.  In particular, every $G_N(q)$ belongs to $\QQ(q)$.

We will show that the two normalized generating series satisfy this equation.

\begin{theorem}\label{thm:common}
Both $\Hseries$ and $\Bseries$ satisfy \eqref{eq:commonequation}.
Consequently, for every $N\ge0$, we have
\begin{equation}\label{eq:commoncoefficients}
 F_N(q)=(q)_N C_{\cc,\le N}(q)
       =q^{-\gamma_N}(q)_N G_N(q).
\end{equation}
\end{theorem}
The proof has two branches, each beginning with a combinatorial
identification.
For the first, let $\Rpaths$ consist of the paths in
$\mathcal D_{aN,bN}$ passing through every $(ka,kb)$, $0\le k\le N$.
Write
\begin{equation}\label{eq:pathseries}
 P_N(q):=\sum_{P\in\Rpaths}q^{\dinv(P)},\qquad
 \Pseries(z;q):=\sum_{N\ge0}\frac{q^{\binom N2}P_N(q)}{(q)_N}z^N,
\end{equation}
where $P_0(q)=1$ and $\dinv(P)$ is the whole-path arm-leg statistic
defined in \eqref{eq:hook}.  Section~\ref{sec:return} expands the Gaussian
factors as filter words and proves
\[
 P_N(q)=q^{\kappa_N}F_N(q),\qquad \Hseries=\Pseries.
\]
Sections~\ref{sec:shuffle}--\ref{sec:Euler} then use the compositional
rational shuffle theorem and a regular operator specialization to show
that $\Pseries$ satisfies \eqref{eq:commonequation}.
The prescribed-return paths define $P_N$, whereas all paths define $A_N$.

For the second branch, Section~\ref{sec:graph} defines a generating
series $\Dseries=\sum_{N\ge 0} D_N(q)z^N$ for collections of $N$ disjoint cycles in a weighted rank graph.  Section~\ref{sec:determinant}
uses a finite determinantal model and Cramer's cofactor expansion to prove that $\Dseries$ satisfies \eqref{eq:commonequation}.
The weight-preserving bijection in Section~\ref{sec:bijection} identifies
collections of $N$ disjoint cycles with balanced cylinders whose entries
are at most $N$.  It gives
\[
 \Dseries(z;q)=\sum_{N\ge0}q^{\gamma_N}C_{\cc,\le N}(q)z^N=\Bseries(z;q).
\]
Therefore, $\Hseries=\Pseries$ and $\Bseries=\Dseries$ both satisfy the common
equation.  Uniqueness identifies the two branches.  
Finally, Section~\ref{sec:conclusion} establishes the main theorems and collects the combinatorial consequences. A finite example and the Lean certificate
are supplied in Sections~\ref{sec:examples} and~\ref{sec:AxiomProver}.

\section*{Acknowledgements} \noindent The authors thank George Andrews and Seewoo Lee for discussions related to this paper. Furthermore, the authors acknowledge Ole Warnaar for sharing his observation that the $q$-series in the HJO Conjecture could be framed in terms of cylindric partitions. 
The authors also thank Eugene Gorsky for suggesting a connection between the coprime high-rank geometry and noncoprime Catalan combinatorics.

\section{Gap coordinates and Gaussian polynomials}\label{sec:definitions}

To interpret the HJO sum combinatorially, we describe its gap coordinates
and factor its generalized multinomial into ordinary Gaussian polynomials.

\subsection{The gap diagram}
The Sylvester coordinates identify the gap set bijectively with
\begin{equation}\label{eq:sylvester}
 G\longleftrightarrow
 \{(x,y)\in\ZZ^2:1\le x<a,\ 1\le y<b,\ ay-bx>0\},
 \qquad g=ay-bx.
\end{equation}
We use the standard two-generator semigroup facts
$|G|=\delta=(a-1)(b-1)/2$ and $\max G=f=ab-a-b$ (see \cite[Section~1.2]{Huang} and \cite[Chapter~2]{RA}).

\begin{figure}[!htbp]
\centering
\begin{tikzpicture}[x=.8cm,y=.8cm,font=\small]
 \draw[step=1,black!20,very thin] (0,0) grid (4,5);
 \foreach \x/\y in {2/4,1/4,1/3}
  {\fill[black!12] (\x-1,\y) rectangle (\x,\y+1);}
 \draw[dashed,black!65] (0,0)--(4,5);
 \foreach \x/\y/\g in {3/4/1,2/4/6,2/3/2,1/4/11,1/3/7,1/2/3}{
  \draw (\x-1,\y) rectangle (\x,\y+1);
  \node at (\x-.5,\y+.5) {$\g$};
 }
 \draw[very thick] (0,3)--(1,3)--(1,4)--(2,4)--(2,5)--(0,5)--cycle;
 \draw[-{Stealth[length=1.8mm]}] (0,0)--(4.7,0) node[right] {$x$};
 \draw[-{Stealth[length=1.8mm]}] (0,0)--(0,5.65) node[above] {$y$};
 \node[below left] at (0,0) {$0$};
 \node[below] at (4,0) {$4$};\node[left] at (0,5) {$5$};
 \fill (1,4) circle (1.9pt);
 \draw[thin] (1,4)--(4.65,3.4);
 \node[anchor=west,align=left] at (4.7,3.4)
    {$(1,4)$: southeast corner\\of the box $11$};
\end{tikzpicture}
\caption{The gap diagram for $(a,b)=(4,5)$ in the above-diagonal
convention.  The shaded filter is $\{6,7,11\}$.
The gap $11=4\cdot4-5\cdot1$ has coordinates $(1,4)$.}
\label{fig:younggap}
\end{figure}
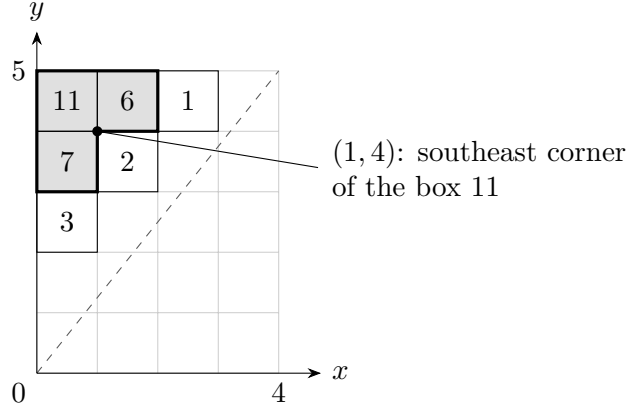

We adopt the convention that the rank of a lattice point $(x,y)$ is $ay-bx$.
A box has the coordinates of its southeast corner, so the box at $[x-1,x]\times[y,y+1]$ has coordinate $(x,y)$ and rank $ay-bx$.
The gap boxes are precisely the full boxes above the diagonal of
$[0,a]\times[0,b]$.  They form a northwest-justified Young diagram, 
as in Figure~\ref{fig:younggap}.
Recall $G$ is a poset by $g\preceq h\iff h-g\in \Gamma$. In Sylvester coordinates, we have
\[
 ay-bx\preceq_G ay'-bx'
       \quad\Longleftrightarrow\quad x\ge x'\text{ and }y\le y'.
\]
Therefore, an order filter is a northwest-justified Young subdiagram.
The maximum gap $f$ occupies the box with southeast corner $(1,b-1)$.
Our horizontal rank increment is $b$, whereas it is $a$ in
\cite[Section~1.2]{Huang}.

\subsection{A polynomial factorization}
We next recall the Gaussian polynomials and use HJO's good-traverse factorization to prove that the generalized multinomial \eqref{eq:mult} is a polynomial with nonnegative coefficients. For integers $u,v$, we recall the Gaussian binomial coefficients
\[
 \qbin uv:=\begin{cases}
\frac{(q)_u}{(q)_v(q)_{u-v}},&0\le v\le u,\\
 0,&\text{otherwise}.
 \end{cases}
\]
\noindent
For $0\le v\le u$, $\qbin uv$ is the inversion generating
polynomial of binary words of length $u$ with $v$ ones.  An {\it inversion}
is a one preceding a zero.  In particular, $\qbin uv\in\NN[q]$.

For the next formula, define the extended upper-boundary value
\begin{equation}\label{eq:boundaryextension}
 \widehat n_j:=\begin{cases}
 n_j,&j\in G,\\ N,&j\in\Gamma,\\0,&j<0.
 \end{cases}
\end{equation}
HJO's good-traverse factorization is
\begin{equation}\label{eq:good}
 \mult N\nn=\prod_{g\in G}
 \qbin{\widehat n_{g+b}-\widehat n_{g-a}}{n_g-\widehat n_{g-a}}.
\end{equation}
This is \cite[Lemma~5.9, (5.20), and Remark~6.15]{HJO}, specialized to the gap poset, using the symmetry $\qbin uv=\qbin u{u-v}$ to write the selected lower dimension.
Visit the gap cells column by column in increasing $x$ and, within each column, increasing $y$. At $g$, the necessary lower predecessor $g-a$ and upper successor $g+b$ have already been visited, or lie on the prescribed boundary. Every earlier lower constraint is contained in the lower predecessor, and every earlier upper constraint contains the upper successor.  This is the ``good traverse'' property used in HJO.

For $\nn\in\mathcal C$ with $n_f\le N$, every Gaussian factor has integer entries satisfying 
$$0\le n_g-\widehat n_{g-a}\le\widehat n_{g+b}-\widehat n_{g-a}.
$$
Therefore, \eqref{eq:mult} is a nonzero polynomial in $\NN[q]$. In Section~\ref{sec:return}, we will use \eqref{eq:good} to give an explicit combinatorial expansion: at each coordinate, the Gaussian counts a choice of labels between two already fixed subsets. The dependence of successive choices is exactly why we can't simply replace the generalized multinomial with an ordinary multinomial.
For $g\in G$, a predecessor $g-a$ cannot belong to $\Gamma$.
Thus $\widehat n_{g-a}=n_{g-a}$ with the zero-extension convention used
in \eqref{eq:PG}. 

\begin{example}[The good traverse for $(a,b)=(4,5)$]
Reading Figure~\ref{fig:younggap} from left to right, and each column
from bottom to top, gives $3,7,11,2,6,1$.
The factorization is
\[
 \mult N\nn=
 \qbin N{n_3}\qbin{N-n_3}{n_7-n_3}
 \qbin{N-n_7}{n_{11}-n_7}\qbin{n_7}{n_2}
 \qbin{n_{11}-n_2}{n_6-n_2}\qbin{n_6}{n_1}.
\]
At each box, its lower predecessor and upper successor are either
already visited or supplied by \eqref{eq:boundaryextension}.
\end{example}

\subsection{Convergence conventions}
All finite identities below are identities of polynomials or rational functions in indeterminates. The notation $[z^N]$ means coefficient extraction in a formal series in $z$. 
The symmetric-function calculation takes place in $\QQ(q)[[z]]$.  The determinant calculation takes place in $\ZZ[[q]][[z]]$.  Both rings embed in $\QQ((q))[[z]]$.
 Here $\QQ((q))$ is the field of Laurent series in $q$ with only finitely many negative powers.

For a fixed vector with $n_f=m$, \eqref{eq:mult} at ambient rank $m$ gives
\[
 P_G(\nn)=\frac{\mult m\nn}{(q)_m}.
\]
Therefore, $P_G(\nn)$ has a power-series expansion with nonnegative coefficients and no negative powers of $q$. The bound \eqref{eq:coercivity} limits $m$ for any fixed coefficient of $q$, and then only finitely many vectors can contribute. This gives a direct formal meaning to \eqref{eq:Z}. We will apply the same principle to the determinant: first bound the data that can contribute to a prescribed coefficient, and only then pass to a limit.

The analytic identity for $|q|<1$ follows as well. HJO proves convergence of its sum in this disk. Alternatively, the Pochhammer quotients are bounded in absolute value by a constant depending on $|q|,a,b$, while the quadratic bound dominates the polynomial growth in the number of vectors. The product in \eqref{eq:HJOconj} converges there because its periodic exponents are bounded. Equality of their formal series therefore also gives equality of the analytic functions.
\subsection{Conventions used in the two constructions}
The proof changes models but not the fixed coprime pair. Table~\ref{tab:conventions} gives the conventions at the points where a change could otherwise alter a formula.
\par
\begin{table}[!htbp]
\centering
\renewcommand{\arraystretch}{1.25}
\begin{tabular}{|p{.25\textwidth}|p{.63\textwidth}|}
\hline
Model & Convention and role\\
\hline
Gap vectors & Use $n_j=0$ off $G$ for gap vectors. For $\widehat n_j$, supply missing values by \eqref{eq:boundaryextension}.\\
\hline
Return paths & Above the diagonal, with a return after every $(a,b)$ block. A filter word is concatenated in reverse order.\\
\hline
Parking functions & Above the diagonal, as in Section~\ref{sec:return}. North labels increase up columns, and $u$ marks area.\\
\hline
Rank graph & North and east edges have increments $+a,-b$.  The east-edge weight uses the ending rank, as in \eqref{eq:rankgraphedges}.\\
\hline
Cylinders & Use the outgoing shifts in \eqref{eq:cylinder}. The parameter $N$ bounds the largest entry.\\
\hline
\end{tabular}
\par\smallskip
\caption{Conventions used in the two constructions, listed by model.}
\label{tab:conventions}
\end{table}
\par

The shifts $\gamma_N$ in \eqref{eq:Hseries} have different origins. On the $\Hseries$ side, hooks between
blocks contribute $\kappa_N$ (Section~\ref{sec:return}), and the creation-operator calculation adds $\binom N2$ to give $\gamma_N$ (Section~\ref{sec:shuffle}).  On the $\Bseries$ side, the cycle construction obtains
$\gamma_N$ by removing a staircase from each residue class (Section~\ref{sec:bijection}).

\section{From poset flags to concatenated rational paths}\label{sec:return}

We express the Gaussian factors as sums over words of order filters and
identify their weights with dinv of concatenated paths.  This proves the
equality $\Hseries=\Pseries$ needed for the first branch of
Theorem~\ref{thm:common}.
We will obtain a combinatorial expansion of the form
\[
 \mult N\nn=\sum_{\boldsymbol D\in\binom{[N]}{\nn}_G}
                         q^{\iota(\boldsymbol D)}.
\]
The finite set $\binom{[N]}{\nn}_G$ and the nonnegative integer statistic
$\iota$ are defined in Section~\ref{subsec:filterwords}.
Since $\mult N\nn$ is the cell-counting polynomial of HJO's poset flag
variety \eqref{eq:poset-flag}, its torus-fixed points and attracting-cell
dimensions suggest the required model.  We use a purely combinatorial
definition in which $\iota$ counts inversions.
Section~\ref{subsec:flaggeometry} verifies the geometric interpretation,
which is not needed in the rest of the proof.

We then identify the combined weight
$Q(\nn(\boldsymbol D))+\iota(\boldsymbol D)$ with the dinv of a single
concatenated path, minus the universal shift $\kappa_N$.  This is Theorem~\ref{thm:worddinv}.

\subsection{Filter words and the Gaussian expansion}\label{subsec:filterwords}
We describe the choices counted by the Gaussian factors and assign their
inversion weights.  Let $\mathcal J(G)$ be the set of order filters of $G$, including the
empty filter.  For $\boldsymbol D=(D_1,\ldots,D_N)\in\mathcal J(G)^N$, put
\begin{equation}\label{eq:wordset}
 \nn(\boldsymbol D):=\sum_{i=1}^N\one_{D_i}.
\end{equation}
Fix $\nn\in\mathcal C$ with $n_f\le N$, and write
$[N]=\{1,\ldots,N\}$, with $[0]=\varnothing$.
We define three equivalent models for the finite set $\binom{[N]}{\nn}_{G}$:
\begin{enumerate}
 \item Words of filters $\{\boldsymbol D\in\mathcal J(G)^N:\nn(\boldsymbol D)=\nn\}$.
 \item Nested subsets $\{(S_g)_{g\in G}:S_g\subseteq[N], |S_g|=n_g, S_g\subseteq S_h\text{ whenever }g\preceq h\}$.
 \item Order filters $\mathcal U\subseteq G\times[N]$ whose slice
       over $g\in G$ has $n_g$ elements, where $[N]$ has the antichain order.
\end{enumerate}
The bijections are
\begin{equation}\label{eq:setslices}
 S_g=\{i:g\in D_i\},\qquad
 D_i=\{g:i\in S_g\},\qquad
 \mathcal U=\{(g,i):g\in D_i\}=\{(g,i):i\in S_g\}.
\end{equation}
The numerical order $1<\cdots<N$ is used for inversions, not to impose
nesting on $D_1,\ldots,D_N$.   

To define $\iota$, use the word model.  For $D,E\in\mathcal J(G)$, set
\begin{equation}\label{eq:I}
 \iota(D,E):=\#\{g\in D:g-a\notin D,\ g\notin E,
                                    \ g+b\in E\cup\Gamma\}.
\end{equation}
More generally, put
\begin{equation}\label{eq:inversionword}
 \iota(\boldsymbol D):=\sum_{1\le i<j\le N}\iota(D_i,D_j).
\end{equation}

The inversion statistic recovers the Gaussian factors for each fixed dimension vector.

\begin{lemma}\label{lem:words}
For $\nn\in\mathcal C$ with $n_f\le N$, we have
\begin{equation}\label{eq:fixedGaussian}
 \mult N\nn=
 \sum_{\substack{\boldsymbol D\in \mathcal{J}(G)^N\\ \nn(\boldsymbol D)=\nn}}
                 q^{\iota(\boldsymbol D)}.
\end{equation}
\end{lemma}
\begin{proof}
Use the good traverse from Section~\ref{sec:definitions}, with boundary
subsets $S_g=[N]$ for $g\in\Gamma$ and $S_g=\varnothing$ for $g<0$.
At gap $g$, the available choice is
$S_{g-a}\subseteq S_g\subseteq S_{g+b}$ with $|S_g|=n_g$.
List $S_{g+b}\setminus S_{g-a}$ in increasing order, marking the selected
labels in $S_g\setminus S_{g-a}$ by one and the others by zero.
The inversion generating polynomial of this binary word is
\[
 \qbin{\widehat n_{g+b}-\widehat n_{g-a}}{n_g-\widehat n_{g-a}}.
\]
Its inversions are exactly the pairs
\[
 i<j,\qquad i\in S_g\setminus S_{g-a},\qquad
 j\in S_{g+b}\setminus S_g.
\]
In the subdiagram word these conditions say
$g\in D_i$, $g-a\notin D_i$, $g\notin D_j$, and
$g+b\in D_j\cup\Gamma$.
Their total is therefore $\iota(\boldsymbol D)$.
Each Gaussian depends on earlier choices only through their prescribed
cardinalities.  Multiplying over the traverse gives \eqref{eq:good},
which proves the identity.
\end{proof}

In the $(4,5)$ example following \eqref{eq:good}, the choice at gap $6$ is
$S_6\setminus S_2$ inside $S_{11}\setminus S_2$.
The corresponding factor is $\qbin{n_{11}-n_2}{n_6-n_2}$.
An inversion $i<j$ satisfies
$6\in D_i$, $2\notin D_i$, $11\in D_j$, and $6\notin D_j$.
This is the contribution of $g=6$ to \eqref{eq:I}.

The full HJO weight is
\begin{equation}\label{eq:wordgrading}
 W(\boldsymbol D):=Q\bigl(\nn(\boldsymbol D)\bigr)
                         +\iota(\boldsymbol D).
\end{equation}
Multiplying \eqref{eq:fixedGaussian} by $q^{Q(\nn)}$ and summing gives
\begin{equation}\label{eq:wordidentity}
 F_N(q)=\sum_{\boldsymbol D\in\mathcal J(G)^N}q^{W(\boldsymbol D)}.
\end{equation}
The same argument retains $t^{|\nn(\boldsymbol D)|}$, or separate
variables for each gap coordinate.  At $q=1$, the lemma gives the
cardinality of the words with any fixed dimension vector.

\subsection{The geometric meaning of inversions}\label{subsec:flaggeometry}
We interpret the inversion statistic as the dimension of an affine cell
in HJO's poset flag variety.  For a field $\Bbbk$ with standard basis $e_1,\ldots,e_N$ of $\Bbbk^N$, and $S\subseteq[N]$, write
$e_S=\operatorname{span}\{e_i:i\in S\}\subseteq\Bbbk^N$. Consider HJO's poset flag variety \cite[Sections~5--6]{HJO}:
\begin{equation}\label{eq:poset-flag}
 \operatorname{Fl}_G(\nn;N)=
 \{(V_g)_{g\in G}:V_g\subseteq\Bbbk^N,\quad
                   \dim V_g=n_g,\quad V_g\subseteq V_h\ (g\preceq h)\}.
\end{equation}
The fixed points of the diagonal torus on $\operatorname{Fl}_G(\nn;N)$ are precisely the coordinate flags  $(e_{S_g})_{g\in G}$, where $(S_g)_{g\in G}$ is a flag of nested subsets in $\binom{[N]}{\nn}_G$.

For a subspace $V$, let $\operatorname{in}(V)$ be the leading-index set
of its reduced column-echelon basis, with the smallest nonzero index
as the pivot.  A subspace with leading set $S$ has basis
\[
 e_i+\sum_{\substack{j>i\\j\notin S}}c_{ji}e_j,\qquad i\in S.
\]
The subspaces with $\operatorname{in}(V)=S$ form the orbit of $e_S$ under the lower triangular Borel $B_-$.  This is the $B_-$-Schubert cell in $\operatorname{Gr}(|S|,N)$. 

At each step of the good traverse, we choose an intermediate subspace
between two previously chosen subspaces.  The following lemma gives
coordinates for this choice when the leading-index sets are prescribed.

\begin{lemma}\label{lem:relativecell}
Let $V_-\subseteq V_+\subseteq\Bbbk^N$ have leading sets $S_-$ and
$S_+$.  For $S_-\subseteq S\subseteq S_+$, the variety
\[
 \{V:V_-\subseteq V\subseteq V_+,\quad\operatorname{in}(V)=S\}
\]
is an affine space of dimension
\[
 \#\{(i,j):i\in S\setminus S_-,\quad j\in S_+\setminus S,\quad i<j\}.
\]
These coordinates vary regularly with the reduced-echelon coordinates
of the bounding pair.
\end{lemma}
\begin{proof}
The Bruhat decomposition for the two-step flag $V_-\subseteq V_+$
reduces the pair to $e_{S_-}\subseteq e_{S_+}$ by a lower triangular
change of basis.  A regular unitriangular representative $U$ is explicit:
its $i$th column is the reduced basis vector of $V_-$ for $i\in S_-$,
that of $V_+$ for $i\in S_+\setminus S_-$, and $e_i$ otherwise.
These columns give $Ue_{S_-}=V_-$ and $Ue_{S_+}=V_+$.
Lower triangular invertible maps preserve leading-index sets.
For the coordinate pair, the desired subspaces are exactly
\[
 e_{S_-}+\operatorname{span}\left\{
 e_i+\sum_{\substack{j\in S_+\setminus S\\j>i}}c_{ji}e_j:
 i\in S\setminus S_-\right\}.
\]
The displayed coefficients are arbitrary and unique.
Applying $U$ gives the assertion.  Both $U^{-1}$ and the ensuing
unit-pivot echelon reductions are polynomial in the input coordinates,
so this parametrization is regular in the bounding pair as well.
\end{proof}

Applying this construction along the traverse gives the geometric interpretation of $\iota$.

\begin{corollary}[Schubert cell dimension]\label{cor:flagcells}
The stratum of $\operatorname{Fl}_G(\nn;N)$ with
$\operatorname{in}(V_g)=S_g$ for every $g$ is an affine space of dimension
$\iota(\boldsymbol D)$.  For any fixed $\lambda_1<\cdots<\lambda_N$, $\lambda_i\in\ZZ$, this stratum is the attracting locus of $(e_{S_g})$ for the one-parameter subgroup
\(
 \operatorname{diag}(\tau^{\lambda_1},\ldots,\tau^{\lambda_N}) \text{ as }\tau\to0.
\)
\end{corollary}
\begin{proof}
Extend the boundary subspaces by $V_g=0$ for $g<0$ and
$V_g=\Bbbk^N$ for $g\in\Gamma$.
At each step of the good traverse, apply Lemma~\ref{lem:relativecell}
to $V_{g-a}\subseteq V_g\subseteq V_{g+b}$.
Its regular parametrization gives a trivial affine-space factor over
the stratum already constructed.  The sum of the dimensions is the
inversion count in Lemma~\ref{lem:words}.
The torus multiplies an echelon coefficient $c_{ji}$ by
$\tau^{\lambda_j-\lambda_i}$, which tends to zero because $j>i$.
Thus its limit is the coordinate flag.
\end{proof}

\subsection{Concatenation and the whole-path statistic}
We now join the paths associated to a filter word and compute the dinv
of their concatenation.  The southeast boundary of a filter $D$ gives a path $P_D$ from $(0,0)$
to $(a,b)$ weakly above the diagonal, with $D$ as its diagram of boxes
above the path.   For a filter word
$\boldsymbol D=(D_1,\ldots,D_N)$, use the reversed concatenation
\begin{equation}\label{eq:upperconcatenation}
 P:=P_{D_N}P_{D_{N-1}}\cdots P_{D_1}.
\end{equation}
It bijects $\mathcal J(G)^N$ with $\Rpaths$, the above-diagonal paths
through every $(ka,kb)$ for $0\le k\le N$.
The reversal is necessary because the inversion statistic retains the
order of the filters.

For an above-diagonal path $P$ of dimensions $(aN,bN)$, let $Y_P$
be its diagram of boxes above the path in the rectangle.
For $\square\in Y_P$, let $\arm_P(\square)$ count the boxes to its right
in the same row, and let $\leg_P(\square)$ count those below it in the
same column.  Neither count includes $\square$ itself.  Define
\begin{equation}\label{eq:hook}
 \dinv(P):=\#\{\square\in Y_P:
 b\arm_P(\square)\le a(\leg_P(\square)+1),\quad
 a\leg_P(\square)<b(\arm_P(\square)+1)\}.
\end{equation}
In other words, $\dinv(P)$ counts boxes whose ``hook slope'' interval $(\leg/(\arm+1),(\leg+1)/\arm]$ contains the diagonal slope $b/a$, where the upper endpoint is treated as $\infty$ if $\arm=0$.
This is the unlabelled rational dinv statistic and the choice of the nonstrict inequality follows the convention of
\cite[Algorithm~6.1(3), p.~21]{BGLX}.
The nonstrict inequality matters in the noncoprime rectangle when $N>1$. In the coprime case $N=1$, endpoint equality is excluded, so this matches the ordinary dinv \cite[Section~1.2]{Huang}, except that arm and leg are interchanged because our horizontal gap increments are $b$ rather than $a$ there. 
Recall $\kappa_N=(a+b-1)\binom N2$.

The combined quadratic and inversion weight has the following
interpretation on the entire concatenation.

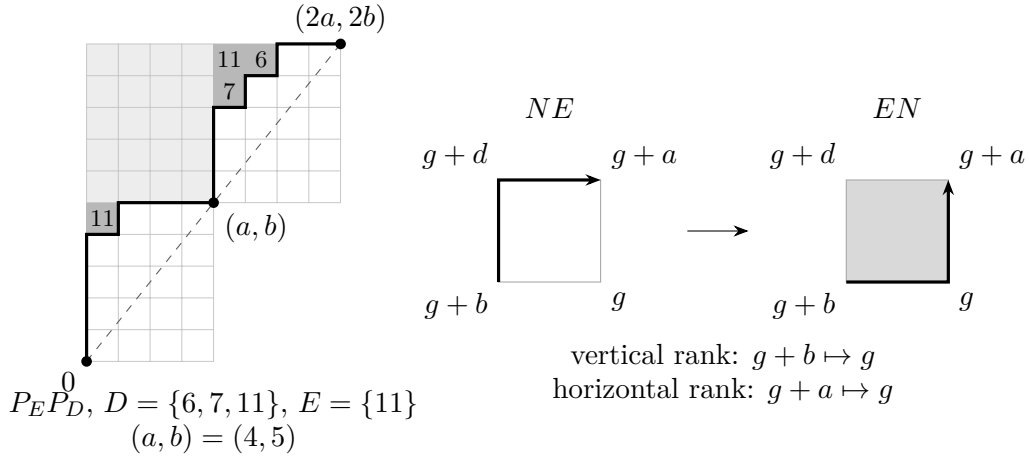
\begin{figure}[!htbp]
\centering
\begin{tikzpicture}[font=\small]
 \begin{scope}[x=.42cm,y=.42cm]
  \fill[black!7] (0,5) rectangle (4,10);
  \draw[step=1,black!22,very thin] (0,0) grid (4,10);
  \draw[step=1,black!22,very thin] (4,5) grid (8,10);
  \draw[dashed,black!65] (0,0)--(8,10);
  \foreach \x/\y/\g in {5/9/6,4/9/11,4/8/7,0/4/11}{
   \fill[black!28] (\x,\y) rectangle (\x+1,\y+1);
   \node[font=\scriptsize] at (\x+.5,\y+.5) {$\g$};
  }
  \draw[very thick] (0,0)--(0,4)--(1,4)--(1,5)--(4,5)--(4,8)
        --(5,8)--(5,9)--(6,9)--(6,10)--(8,10);
  \fill (0,0) circle (2pt) (4,5) circle (2pt) (8,10) circle (2pt);
  \node[below left] at (0,0) {$0$};
  \node[below right] at (4,5) {$(a,b)$};
  \node[above] at (8,10) {$(2a,2b)$};
  \node[below,align=center] at (4,-.55)
    {$P_EP_D$, $D=\{6,7,11\}$, $E=\{11\}$\\$(a,b)=(4,5)$};
 \end{scope}
 \begin{scope}[xshift=5.45cm,yshift=1.05cm,x=1.35cm,y=1.35cm]
  \draw[black!35] (0,0) rectangle (1,1);
  \draw[very thick,-{Stealth[length=2mm]}] (0,0)--(0,1)--(1,1);
  \node[above left] at (0,1) {$g+d$};
  \node[above right] at (1,1) {$g+a$};
  \node[below left] at (0,0) {$g+b$};
  \node[below right] at (1,0) {$g$};
  \node at (.5,1.7) {$NE$};
  \draw[-{Stealth[length=2mm]}] (1.85,.5)--(2.45,.5);
  \begin{scope}[xshift=4.6cm]
   \fill[black!15] (0,0) rectangle (1,1);
   \draw[black!35] (0,0) rectangle (1,1);
   \draw[very thick,-{Stealth[length=2mm]}] (0,0)--(1,0)--(1,1);
   \node[above left] at (0,1) {$g+d$};
   \node[above right] at (1,1) {$g+a$};
   \node[below left] at (0,0) {$g+b$};
   \node[below right] at (1,0) {$g$};
   \node at (.5,1.7) {$EN$};
  \end{scope}
  \node[align=center] at (2.2,-.9)
   {vertical rank: $g+b\mapsto g$\\horizontal rank: $g+a\mapsto g$};
 \end{scope}
\end{tikzpicture}
\caption{The reversed concatenation for the word $(D,E)$.
The light rectangle contains its $ab$ cross-block boxes, and darker
boxes are internal to the two blocks.  Adding a box of southeast rank
$g$ changes $NE$ to $EN$ and gives the boundary-rank changes on the right.}
\label{fig:gap}
\end{figure}

\begin{theorem}\label{thm:worddinv}
For $\boldsymbol D=(D_1,\ldots,D_N)\in\mathcal J(G)^N$ and
$P=P_{D_N}\cdots P_{D_1}$ as in \eqref{eq:upperconcatenation}, we have
\begin{align}
 \dinv(P)&=\kappa_N+W(\boldsymbol D),\label{eq:pointwisedinv}\\
 \area(P)&=N\delta-\sum_{i=1}^N|D_i|.\label{eq:wordarea}
\end{align}
\end{theorem}
The area identity follows by addition over the blocks, each of which has
area $\delta-|D_i|$.  For the dinv identity, we must also count hooks
between different blocks.  Figure~\ref{fig:gap} shows the concatenation
and the local boundary change used in the proof.

\begin{proof}
The rank of a horizontal step is the rank of its east endpoint.
The rank of a vertical step is the rank of its south endpoint.
For a box $\square$ with southeast coordinates $(x,y)$, the vertical
step at the end of its arm starts at $(x+\arm_P(\square),y)$.
The horizontal step at the end of its leg ends at
$(x,y-\leg_P(\square))$.
Write $v(\square)$ and $h(\square)$ for these two step ranks.
Since the rank is $ay-bx$,
\begin{equation}\label{eq:rankhook}
 \begin{aligned}
 v(\square)-h(\square)
     &=a\leg_P(\square)-b\arm_P(\square),\\
 \square\text{ is counted in \eqref{eq:hook}}
     &\quad\Longleftrightarrow\quad -a\le v(\square)-h(\square)<b.
 \end{aligned}
\end{equation}
The horizontal step precedes the vertical step.
Conversely, every such ordered pair determines one box.
Figure~\ref{fig:hookranks} records this correspondence.

\begin{figure}[!htbp]
\centering
\begin{tikzpicture}[x=.68cm,y=.68cm,font=\small]
 \fill[black!10] (2,3) rectangle (5,4);
 \fill[black!10] (1,0) rectangle (2,3);
 \fill[black!27] (1,3) rectangle (2,4);
 \draw[step=1,black!25,very thin] (1,3) grid (5,4);
 \draw[step=1,black!25,very thin] (1,0) grid (2,3);
 \draw[very thick] (0,0)--(2,0)--(2,3)--(5,3)--(5,5)--(7,5);
 \draw[very thick,-{Stealth[length=1.7mm]}] (5,3)--(5,4);
 \draw[very thick,-{Stealth[length=1.7mm]}] (1,0)--(2,0);
 \fill (2,0) circle (1.7pt) (2,3) circle (1.7pt) (5,3) circle (1.7pt);
 \node at (1.5,3.5) {$\square$};
 \node[above] at (3.5,4) {$\arm_P(\square)$};
 \node[left] at (1,1.5) {$\leg_P(\square)$};
 \node[anchor=west,align=left] at (5.3,2.5)
    {vertical rank $v$\\$(x+\arm_P(\square),y)$};
 \node[below,align=center] at (2,-.3)
    {horizontal rank $h$\\$(x,y-\leg_P(\square))$};
 \node[anchor=west] at (2.2,2.6) {$(x,y)$};
\end{tikzpicture}
\caption{The hook of an above-path box.  Its arm extends rightward and
its leg downward.  The southeast corner is $(x,y)$.
The horizontal step precedes the vertical step.}
\label{fig:hookranks}
\end{figure}
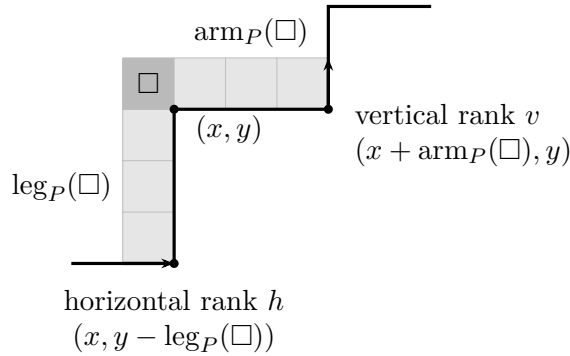

We encode boundary ranks by Laurent series, using the technique of \cite[Lemma~6.16]{HJO}.  Its basic generating identity is
\begin{equation}\label{eq:gammagenerating}
 \Gamma(T):=\sum_{g\in\Gamma}T^g
          =\frac{1-T^{ab}}{(1-T^a)(1-T^b)},
\end{equation}
see for instance \cite{BrownShiue}.
For completeness, every semigroup element has a unique representation
$au+bv$ with $u\ge0$ and $0\le v<a$.  Summing these monomials proves
the formula.
For this proof write
$D(T)=\sum_{g\in D}T^g$ and $E(T)=\sum_{g\in E}T^g$.
Use the coefficient pairing
\begin{equation}\label{eq:laurentpairing}
 \langle T^r,T^s\rangle=\one_{r=s},\qquad
 \langle T^r A,B\rangle=\langle A,T^{-r}B\rangle.
\end{equation}
At least one argument in every pairing below is a Laurent polynomial.  Therefore, each pairing uses only finitely many coefficients. For any Laurent polynomial $A(T)$, the adjoint of multiplication by $A(T)$ is multiplication by
$A(T^{-1})$. 

Let $\mathcal V_D$ and $\mathcal H_D$ be the multisets of vertical
south-endpoint and horizontal east-endpoint ranks of $P_D$.  Their generating polynomials are
\[
 V_D(T):=\sum_{v\in\mathcal V_D}T^v,\qquad
 H_D(T):=\sum_{h\in\mathcal H_D}T^h.
\]
We suppress $(T)$ in $V_D,H_D$ when there is no ambiguity.
The path of the empty filter has
\[
 V_0=\frac{1-T^{ab}}{1-T^a},\qquad
 H_0=\frac{1-T^{ab}}{1-T^b}.
\]
Adding a box of rank $g$ changes $NE$ to $EN$, as in
Figure~\ref{fig:gap}.  Its vertical rank changes from $g+b$ to $g$,
and its horizontal rank from $g+a$ to $g$.  Hence
\begin{equation}\label{eq:rankmeasures}
 V_D=V_0+(1-T^b)D(T),\qquad
 H_E=H_0+(1-T^a)E(T).
\end{equation}
By \eqref{eq:rankhook}, the cross-block contribution to $\dinv(P_EP_D)$ is
\begin{equation}\label{eq:crossdefinition}
 C(D,E)=\sum_{v\in\mathcal V_D,\,h\in\mathcal H_E} \boldsymbol{1}_{-a\leq v-h<b} = \langle V_D,LH_E\rangle,\qquad
 L=\sum_{r=-a}^{b-1}T^r=\frac{T^{-a}-T^b}{1-T}.
\end{equation}
The same count applies to each ordered pair of blocks, even when they
are not adjacent.  Therefore, we have
\begin{equation}\label{eq:dinv-in-cross}
    \dinv(P_{D_N}\cdots P_{D_1})=\sum_{i=1}^N \dinv(P_{D_i}) + \sum_{i<j} C(D_i,D_j).
\end{equation}
Indeed, for $i<j$, the block $P_{D_j}$ occurs before $P_{D_i}$ in
\eqref{eq:upperconcatenation}.  Every horizontal step of $P_{D_j}$
therefore precedes every vertical step of $P_{D_i}$.
Translation by $(ka,kb)$ preserves the rank $ay-bx$, so the contribution
is $C(D_i,D_j)$ regardless of the distance between the blocks.

Expanding \eqref{eq:rankmeasures} expresses $C(D,E)$ as a constant,
two linear terms, and a bilinear term.  We compute them in this order.

\medskip
\noindent
\emph{The constant term.}
Since $H_0=(1-T^a)\Gamma(T)$, adjointness gives
$\langle V_0,LH_0\rangle=\langle(1-T^{-a})L(T^{-1})V_0,\Gamma(T)\rangle$.
Cancelling the geometric factors yields
\[ \langle V_0, LH_0\rangle = \langle (T^{ab}-1)\frac{1-T^{-(a+b)}}{1-T^{-1}},\Gamma(T)\rangle,\]
where the right-hand side is a pairing of a Laurent polynomial with a series.
Expanding the geometric sum, the terms of the first argument with
nonnegative exponents are
\[ -1 + (T^{ab-a-b+1}+\dots+T^{ab}).\]
Since $\max G=f=ab-a-b$, all these terms contribute to the pairing with $\Gamma(T)$, so $\langle V_0,LH_0\rangle=a+b-1$.
\medskip

\noindent
\emph{The linear terms.}
Adjointness and cancellation of the finite geometric sums give
\begin{align*}
 \langle(1-T^b)D,LH_0\rangle
   &=\left\langle D,(T^{ab}-1)\sum_{k=1}^{d}T^{-k}\right\rangle,\\
 \langle V_0,L(1-T^a)E\rangle
   &=\left\langle (T^{ab}-1)\sum_{k=0}^{d-1}T^{-k},E\right\rangle.
\end{align*}
The first Laurent polynomial has exactly one term supported on $[1,f]$,
namely $T^f$.  The second has no term supported on $[1,f]$.
Since $D,E$ are supported in $[1,f]$, the linear contributions are
$\one_D(f)$ and $0$.
\medskip

\noindent
\emph{The bilinear term.}
We define the polynomial
\[
 K(T):=\sum_{s\in\ZZ}K(s)T^s=\frac{(1-T^a)(1-T^b)}{1-T}.
\]
Adjointness gives
\[ \langle (1-T^b)D(T),L(1-T^a)E(T)\rangle = \langle D(T), K(T)(1-T^{-d})E(T)\rangle.\]
Combining the constant, linear, and bilinear terms gives
\begin{equation}\label{eq:cross-simplified}
    C(D,E)=a+b-1+\one_D(f)+\langle D(T), K(T)(1-T^{-d})E(T)\rangle.
\end{equation}

We now compare this expression with $W(\boldsymbol D)$.
By bilinear expansion and the identity $Q(1_D)=\dinv(D)$ from \cite[Theorem~1.1]{Huang}, we have
\[ Q(\one_{D_1}+\dots+\one_{D_N})=\sum_{i=1}^N \dinv(P_{D_i})+\sum_{i<j} 2\mathcal{B}(D_i,D_j),\]
where
$$2\mathcal{B}(D,E)\coloneqq \sum_{g,h\in G} (K(g-h)+K(h-g))\,\one_{D}(g)\one_{E}(h).$$
Recall \eqref{eq:inversionword}, \eqref{eq:wordgrading}, and \eqref{eq:dinv-in-cross}. By cancelling the contribution of $\sum_{i=1}^N \dinv(P_{D_i})$ and $\kappa_N=\binom{N}{2}(a+b-1)$, the goal identity \eqref{eq:pointwisedinv} reduces to a pairwise identity
\begin{equation}\label{eq:pairwise-certificate}
    \one_D(f)+\langle D(T), K(T)(1-T^{-d})E(T)\rangle=2\mathcal{B}(D,E)+\iota(D,E).
\end{equation}
To prove this identity, first note that
\begin{align*}
    2\mathcal{B}(D,E)&=\langle D(T),(K(T)+K(T^{-1}))E(T)\rangle.
\end{align*}
On the other hand, the definition gives
\begin{align*}
    \iota(D,E)&=\sum_{g\in G} (\one_D(g)-\one_D(g-a))(\widehat\one_E(g+b)-\widehat\one_E(g)), \qquad \widehat\one_E=\one_E+\one_\Gamma.
\end{align*}
Observe that the sum can be extended to $\ZZ$ without change: if $g<0$, then the first factor vanishes, and if $g\in \Gamma$, then the second factor vanishes. As a result, we have
\begin{align*}
    \iota(D,E)&=\sum_{g\in \ZZ} (\one_D(g)-\one_D(g-a))(\widehat\one_E(g+b)-\widehat\one_E(g))\\
    &=\langle (1-T^a)D(T),(T^{-b}-1)(E(T)+\Gamma(T))\rangle.
\end{align*}
Polynomial arithmetic with adjointness gives
\begin{align*}
    \iota(D,E)&=\langle D(T),T^f-T^{-(a+b)}\rangle-\langle D(T),(1-T^{-a})(1-T^{-b})E(T)\rangle \\
    &=\one_D(f)-\langle D(T),(1-T^{-a})(1-T^{-b})E(T)\rangle,
\end{align*}
where the last equality uses $D\subseteq G$.
Finally, we  have
\[
 K(T^{-1})=-T^{1-d}K(T),\qquad
 (1-T^{-a})(1-T^{-b})=T^{-d}(1-T)K(T).
\]
Therefore, we have
\[
 K(T)+K(T^{-1})-(1-T^{-a})(1-T^{-b})
       =K(T)(1-T^{-d}).
\]
Substitution proves \eqref{eq:pairwise-certificate}.
Combining this with \eqref{eq:dinv-in-cross} proves
\eqref{eq:pointwisedinv}.
\end{proof}

Summing the pointwise identity gives the return-path model, with the
area statistic retaining the variable $t$.

\begin{corollary}\label{prop:return}
For every $N\ge0$, we have
\begin{align}
 q^{\kappa_N}F_N(q)&=\sum_{P\in\Rpaths}q^{\dinv(P)},\label{eq:return}\\
 N_{a,b;N}(q,t)&=q^{-\kappa_N}t^{N\delta}
          \sum_{P\in\Rpaths}q^{\dinv(P)}t^{-\area(P)}.\label{eq:returnfull}
\end{align}
\end{corollary}
\begin{proof}
Multiply the fixed-vector identity \eqref{eq:fixedGaussian} by
$q^{Q(\nn)}t^{|\nn|}$, sum over $\nn$, and apply
Theorem~\ref{thm:worddinv}.  Setting $t=1$ gives \eqref{eq:return}.
\end{proof}
In the notation of \eqref{eq:pathseries}, the first identity is
$P_N=q^{\kappa_N}F_N$.  Since
$\gamma_N=\kappa_N+\binom N2$, it gives
\begin{equation}\label{eq:HP}
 \Hseries(z;q)=\Pseries(z;q).
\end{equation}

\section{A scalar homomorphism from symmetric functions}\label{sec:shuffle}
To prove the $q$-difference equation for $\Pseries$, we construct an
algebra homomorphism $\Phi$ from symmetric functions to $\QQ(q)$.
In Section~\ref{sec:Euler}, its values on two symmetric functions will
give the polynomials $P_N(q)$ and $A_N(q)$.

We use the slope operators of Bergeron--Garsia--Leven--Xin
\cite{BGLX}.  Their values are initially operators over $\QQ(q,u)$.
Evaluation on $1$ becomes multiplicative only after a regular
specialization at $u=1$, which we prove below.
The source locators for BGLX and Mellit refer to arXiv versions~2 and~1,
respectively.  Symmetric-function conventions follow
\cite[Chapter~I, Sections~2--4]{Macdonald}.

\subsection{Symmetric functions and the slope homomorphism}
We recall the plethystic conventions and BGLX operators needed to define the slope homomorphism $\Theta_{a,b}$.
 Let $X:=x_1+x_2+\cdots$ denote a formal alphabet, and set $p_k[X]:=\sum_i x_i^k$. The ring of symmetric functions over our coefficient field is
\[
 \Lambda:=\QQ(q,u)[p_1,p_2,\ldots],\qquad \deg p_k:=k.
\]
Each of its elements has bounded degree and hence is a polynomial in finitely many power sums. Define the complete and elementary symmetric functions by
\begin{align}
 H_X(v):=\sum_{n\ge0}h_n[X]v^n
 &=\exp\left(\sum_{k\ge1}\frac{p_k[X]v^k}{k}\right),\label{eq:Hgenerating}\\
 E_X(v):=\sum_{n\ge0}e_n[X]v^n
 &=\exp\left(\sum_{k\ge1}\frac{(-1)^{k-1}p_k[X]v^k}{k}\right).
 \label{eq:Egenerating}
\end{align}
Set $h_n=e_n=0$ for $n<0$ and $h_0=e_0=1$. The Schur basis is given by $s_\lambda:=\det(h_{\lambda_i-i+j})$, with determinant size at least the number of nonzero parts of $\lambda$.
 In particular $s_{(n)}=h_n$ and $s_{(1^n)}=e_n$.

The variables $q,u$ will be the dinv and area parameters.
The HJO variable $t$ in \eqref{eq:returnfull} is separate.
For plethystic substitution, use the alphabet algebra
\[
 R=\QQ(q,u)[p_1,p_2,\ldots][z,z^{-1}].
\]
For each $j\ge1$, define the $\QQ$-algebra homomorphism
$\psi_j:R\to R$ by
\[
 \psi_j(q)=q^j,\quad \psi_j(u)=u^j,\quad
 \psi_j(z)=z^j,\quad \psi_j(p_k)=p_{jk}.
\]
We write $p_j[E]=\psi_j(E)$ for $E\in R$, identifying the alphabet
$X$ with $p_1$.  Thus $p_j[X]=p_j$, whereas $p_j[2]=2$,
$p_j[q]=q^j$, and $p_j[-E]=-p_j[E]$.  The maps $\psi_j$ fix $\QQ$, not $\QQ(q,u)$.

If $f=\sum_\lambda c_\lambda(q,u) p_\lambda\in\Lambda$, then its value on $E$ is
\[
 f[E]=\sum_\lambda c_\lambda\prod_i p_{\lambda_i}[E].
\]
We let $M=(1-q)(1-u).$
These conventions give
\[
 p_j[X/(1-q)]=\frac{p_j}{1-q^j},\qquad
 p_j[M/z]=(1-q^j)(1-u^j)z^{-j}.
\]
For example, we have
\[
 p_2[X-(1-q^{-1})/z]=p_2-(1-q^{-2})z^{-2}.
\]
These are the standard plethystic conventions
\cite[Chapter~1, Plethystic Notation]{Haglund}.
Each alphabet substitution below is a finite expression in $R$.
Multiplication by $H_X(z)$ or $E_X(-z)$ gives an ordinary formal
Laurent series, in which we extract coefficients.
The variable $v$ in \eqref{eq:Hgenerating}--\eqref{eq:Egenerating}
is an ordinary series variable, not part of a plethystic argument.

We next introduce the operators that allow the shuffle theorem to be applied. Their recursion is also needed later to prove regularity at $u=1$. BGLX's basic operators are
\begin{equation}\label{eq:Dbase}
 \Dop_n f[X]:=[z^n]\left(f[X+M/z]E_X(-z)\right), \quad n\geq 0.
\end{equation}
The operator $\Dop_0$ is already complicated; see \eqref{eq:firstorderD}. An operator product means composition, and $[P,Q]=PQ-QP$. Coefficient extraction in \eqref{eq:Dbase} is finite on each input: the plethysm is Laurent polynomial in $z$, making the series before coefficient extraction a Laurent series of $z$ with coefficients in $\Lambda$. The operator $\Dop_n$ raises homogeneous degree by $n$.

We now define the operators needed in the proof. They are $\Qop_{m,n}$, $(m,n)\in (\NN\times \ZZ_{\ge 1})\cup \{(1,0)\}$. We start with the leftmost two columns
\[
 \Qop_{1,n}:=\Dop_n\quad(n\ge0),\qquad
 \Qop_{0,n}f:=U_n f\quad(n\ge1),
\]
where
\begin{equation}\label{eq:Uk}
 U_n[X;q,u]:=\frac{qu}{qu-1}
        h_n\left[\frac{1-qu}{qu}X\right].
\end{equation}
Therefore, $\Qop_{0,n}$ is multiplication by $U_n$, starting with $U_1=-e_1$ and $U_2=((qu-1)p_1^2-(qu+1)p_2)/(2qu)$. This is the vertical-axis
normalization of \cite[Algorithm~3.1 and Remark~5.1]{BGLX}.

For $m,n\geq 1$, we say $(m,n)=(m'',n'')+(m',n')$ is a \emph{split}, written $\operatorname{Split}(m,n)=(m'',n'')+(m',n')$, if
\[ m=m'+m'',\ n=n'+n'',\ m',n',m'',n''\ge 0, \quad \det\begin{bmatrix} m & n \\ m' & n'\end{bmatrix}=mn'-nm'=\gcd(m,n).\]
It follows that $\gcd(m',n')=\gcd(m'',n'')=1$. In other words, we pick a splitting point $(m'',n'')$ strictly below the diagonal from $(0,0)$ to $(m,n)$ that is as close as possible to the diagonal. 
If $m,n\ge 1$ are coprime, then the split is unique. For example,
\[
 \operatorname{Split}(m,1)=(1,0)+(m-1,1)\quad(m\ge 1), \qquad \operatorname{Split}(1,n)=(1,n-1)+(0,1)\quad(n\ge 1).
\]
If $m,n>1$, then $0<m',m''<m$ and $0<n',n''<n$. If $\gcd(m,n)=k>1$, then the split is not unique; we can pick the splitting point from any coprime block, for example, $(6,10)$ has the two splits $(2,3)+(4,7)$ and $(5,8)+(1,2)$. Concretely, for $\gcd(m,n)=1$, the splits of $(km,kn)$ are
\[ \operatorname{Split}_j(km,kn)=((j-1)m+m'',(j-1)n+n'')+((k-j)m+m',(k-j)n+n'), \qquad 1\le j\le k.\]

Definition~2.1 and Theorem~5.1 of \cite{BGLX} show $\Qop_{m,n}, m,n\ge 1$ are well defined recursively by
\begin{equation}\label{eq:Qrecursion}
 \Qop_{m,n}:=\frac{[\Qop_{m',n'},\Qop_{m'',n''}]}{M}
       , \qquad m,n\ge 1,\ \operatorname{Split}(m,n)=(m'',n'')+(m',n')
\end{equation}
with boundary cases $\Qop_{0,1}$ and $\Qop_{1,0}$ already defined. The content of the well-definedness is the following. For coprime $m,n\ge 1$, the split is unique. The boundary case gives a recursive algorithm to define $\Qop_{m,n}$ for coprime $m,n\geq 1$. This recursively defines $\Qop_{1,n}$, and its value is $\Qop_{1,n}=\Dop_n$ by \cite[(4.9)(ii)]{BGLX}.\footnote{BGLX take $\Qop_{1,n}=\Dop_n$ as a boundary definition in \cite[(2.5) and Definition~2.1]{BGLX}, rather than defining it recursively. Since $\operatorname{Split}(1,n)=(1,n-1)+(0,1)$ and $\Qop_{0,1}=-e_1$, the consistency between the boundary definition and the recursive definition is equivalent to $\frac{[\Dop_{n-1},e_1]}{M}=\Dop_n$, which is \cite[(4.9)(ii)]{BGLX}.} For non-coprime $m,n\ge 1$, Theorem~5.1 of \cite{BGLX} shows that the right-hand side of \eqref{eq:Qrecursion} is independent of the choice of the split. Thus, we may define $\Qop_{m,n}$ using any split and the already defined coprime operators. Theorem~5.1 of \cite{BGLX} also shows collinear operators commute, namely, we have
\[
 [\Qop_{km,kn},\Qop_{k'm,k'n}]=0\qquad k,k',m,n\ge1.
\]
From the definitions above, all operators $\Qop_{m,n}$, $(m,n)\in (\NN\times \ZZ_{\ge 1})\cup \{(1,0)\}$ are generated by scalar multiplications by $U_n$, the operator $\Dop_0$, and the operation $(P,Q)\mapsto [P,Q]/M$. The operator $\Qop_{m,n}$ raises homogeneous degree by $n$.

The symmetric functions $U_1,U_2,\dots$ freely generate
$\Lambda$ as a $\QQ(q,u)$-algebra; see Lemma~\ref{lem:u-p-transition}.  Collinear commutativity gives the $\QQ(q,u)$-algebra homomorphism
\[
 \Theta_{a,b}:\Lambda\longrightarrow
       \operatorname{End}_{\QQ(q,u)}(\Lambda),\qquad
 \Theta_{a,b}(U_n):=\Qop_{na,nb}.
\]
This is the construction in \cite[Algorithm~3.1]{BGLX}.
We write $\Theta=\Theta_{a,b}$ when the pair is understood.
The map fixes $\QQ(q,u)$, and $\Theta(f)$ is a
$\QQ(q,u)$-linear operator on symmetric functions.  If $f$ is
homogeneous of degree $N$, then $\Theta(f)$ raises degree by $bN$.

\subsection{A differential-operator argument at \texorpdfstring{$u=1$}{u=1}}\label{sec:rees}
We will show that for a symmetric function $f\in \QQ(q)[p_1,p_2,\dots]$, the operator $\Theta(f)$ is regular at $u=1$ and specializes to a
multiplication operator (which forces the multiplier to be its value on $1$). We will interpret any operator we encounter as a differential operator in the $p_j$ variables, so it is a multiplication operator if and only if its differential order is zero. The difficulty is that the recursion \eqref{eq:Qrecursion} requires division
$M=(1-q)(1-u)$, which vanishes at $u=1$.

Put $\hbar:=1-u$.  We embed $\QQ(q,u)$ into $\QQ(q)[[\hbar]][\hbar^{-1}]$; an element of $\QQ(q,u)$ is regular if it is in $\QQ(q)[[\hbar]]$, i.e., its Laurent expansion in $\hbar$ has no negative powers.  It would suffice at
the end to know that the constant term of each operator has
differential order zero.  To prove this through the recursion, we
keep track of the order of every coefficient.

The derivatives arise from the ordinary Taylor formula in the power
sums. Let 
\[
 \Delta_j:=p_j[M/z]/\hbar=(1-q^j)\bigl(1-(1-\hbar)^j\bigr)z^{-j}/\hbar
      =j(1-q^j)z^{-j}+O(\hbar);
\]
in particular, it is regular.
For $f=f(p_1,p_2,\dots)\in \QQ(q)[p_1,p_2,\dots]$, the substitution $f[X+M/z]$ is given by
\[ f[X+M/z] = f(p_1[X+M/z],p_2[X+M/z],\dots)=f(p_1+\hbar\Delta_1,p_2+\hbar\Delta_2,\dots)\]
By the Taylor expansion for the polynomial $f(p_1+\hbar\Delta_1,p_2+\hbar\Delta_2,\dots)$ in $p_1,p_2,\dots$, we have
\begin{equation}\label{eq:plethTaylor}
 f[X+M/z]=\sum_\beta\frac{\Delta^\beta}{\beta!}D^\beta f,
 \qquad
 \Delta^\beta:=\prod_{j\ge1}\Delta_j^{\beta_j},\quad D^\beta=\prod_j(\hbar\frac{\partial}{\partial p_j})^{\beta_j},\quad
 \beta!:=\prod_{j\ge1}\beta_j!.
\end{equation}
Here the sum ranges over $\beta$ with finite support. We note that the sum is finite for each $f$. 

Consider the Rees algebra $\mathscr R$ generated over
$\QQ(q)[[\hbar]]$ by multiplication by the $p_j$ and the operators
$\hbar\partial/\partial p_j$, allowing locally finite sums and
coefficientwise formal expansions in $\hbar$. It acts on
$\QQ(q)[p_1,p_2,\dots][[\hbar]]$. At each power of $\hbar$, only
finitely many differential monomials act nontrivially on a fixed
polynomial input. With the brackets understood in this completed
sense, we write
\begin{equation}\label{def:rees}
    \mathscr R=\QQ(q)[[\hbar]][p_1,p_2,\dots][\hbar \partial/\partial_{p_1},\hbar \partial/\partial_{p_2},\dots].
\end{equation}

Recall the definition of $\Dop_n$ from \eqref{eq:Dbase}. Since $E_X(-z)$ is regular in $\hbar$, and the $D^\beta$-coefficient of the operator $f\mapsto f[X+M/z]$ is regular by \eqref{eq:plethTaylor}, we conclude that $\Dop_n\in \mathscr R$.

\begin{remark}
The ring $\mathscr{R}$ can be characterized by an order bound. Write $P=\sum_{r\geq0}\hbar^rP_r$, where each $P_r$ is a polynomial-coefficient differential operator in the variables $p_1,p_2,\dots$. Then $P\in\mathscr R$ if and only if $P_r$ has differential order at most $r$ for every $r$. For example, truncating \eqref{eq:plethTaylor} up to the $\hbar^1$ coefficients gives
\[
 f[X+M/z]=f+\hbar\sum_{j\ge1}j(1-q^j)z^{-j}
               \frac{\partial f}{\partial p_j}+O(\hbar^2).
\]
Since $E_X(z)=\sum_{j\ge 0}e_j z^j$, definition \eqref{eq:Dbase} now gives
\begin{equation}\label{eq:firstorderD}
 \Dop_n=(-1)^n e_n+
 \hbar\sum_{j\ge1}j(1-q^j)(-1)^{n+j}e_{n+j}
              \frac{\partial}{\partial p_j}+O(\hbar^2).
\end{equation}
The constant term has no differential order, and the $\hbar^1$ term only has the first differential, as is expected from the fact $\Dop_n\in \mathscr{R}$.
\end{remark}

Recall that the definition \eqref{eq:Qrecursion} for $\Qop_{m,n}$ requires the operation $(P,Q)\mapsto [P,Q]/M$. 

\begin{lemma}\label{lem:rees}
The ring $\mathscr R$ is closed under the operation
\[
 (P,Q)\longmapsto\frac{[P,Q]}{M}.
\]
Every element of $\mathscr R$ specializes at $u=1$ to multiplication
by its value on the seed $1$.
\end{lemma}
\begin{proof}
The generators commute modulo $\hbar\mathscr R$, since
$[\hbar\partial/\partial p_i,p_j]=\hbar\delta_{ij}$ and the remaining
pairs commute. For the allowed sums and formal expansions, the
product rule gives the same conclusion: the coefficient of $\hbar^r$
in $[P,Q]$ has differential order at most $r-1$ for $r\ge1$, and its
constant coefficient is zero. In particular, for any $P,Q\in \mathscr R$, $[P,Q]\in \hbar \mathscr R$. Since $M=(1-q)\hbar$ and $1-q$ is invertible in $\QQ(q)[[\hbar]]$, this proves the first assertion. In addition, for $P\in \mathscr R$, the $\hbar=0$ specialization of its action $P|_{u=1}$ is a multiplication operator by the order bound, say $P|_{u=1} g=fg$. Testing with $g=1$ gives $f=P|_{u=1}1$, proving the second assertion.
\end{proof}

The closure statement now propagates regularity from the boundary operators through the recursive slope construction.

\begin{lemma}\label{lem:qopregular}
    Every operator $\Qop_{m,n}$ with $(m,n)\in (\NN\times \ZZ_{\ge 1})\cup \{(1,0)\}$ belongs to $\mathscr R$. In particular, we have $\Theta_{a,b}(U_n)\in \mathscr{R}$ for $n\ge 1$.
\end{lemma}
\begin{proof}
    For the boundary cases, we have argued that $\Qop_{1,n}=\Dop_n\in \mathscr R$. Note that $U_n$ is regular at $u=1$: a direct substitution gives
    \begin{equation}\label{eq:UkSpecial}
        U_n[X;q,1]:=\frac{q}{q-1}
        h_n\left[\frac{1-q}{q}X\right]\in \QQ(q)[p_1,p_2,\dots].
    \end{equation}
    Therefore, we have $\Qop_{0,n}\in \mathscr R$. Because any $\Qop_{m,n}$ is generated from these boundary cases via the operation $(P,Q)\mapsto [P,Q]/M$, it is in $\mathscr R$ by Lemma~\ref{lem:rees}.
\end{proof}

We also need that
expressing a symmetric function in the generators $U_n$ doesn't introduce a
pole.  
\begin{lemma}\label{lem:u-p-transition}
    As algebras, we have $\QQ(q)[[\hbar]][p_1,p_2,\dots]=\QQ(q)[[\hbar]][U_1,U_2,\dots]$.
\end{lemma}
\begin{proof}
    By \eqref{eq:Hgenerating}, we have
    \begin{align*}
        \frac{qu}{qu-1}+\sum_{n\geq 1} U_n[X]v^n&=\frac{qu}{qu-1}\sum_{n\geq 0} h_n\left[\frac{1-qu}{qu}X\right]v^n=\frac{qu}{qu-1}\exp\left( \sum_{k\geq 1} \frac{p_k\left[\frac{1-qu}{qu}X\right]}{k}v^k\right)\\
        &=\frac{qu}{qu-1}\exp\left( \sum_{k\geq 1} \frac{(q^{-k}u^{-k}-1)p_k\left[X\right]}{k}v^k\right),
    \end{align*}
    and inversely, we have
    \begin{align*}
        \sum_{k\geq 1} \frac{(q^{-k}u^{-k}-1)p_k\left[X\right]}{k}v^k&=\log\left(1+\frac{qu-1}{qu}\sum_{n=1}^\infty U_n[X]v^n\right).
    \end{align*}
    All the coefficients involved in the two equations above are regular and nonvanishing at $u=1$. Hence, direct expansion gives explicit transition polynomials, showing $p_k\in \QQ(q)[[\hbar]][U_1,\dots,U_k]$ and $U_n\in \QQ(q)[[\hbar]][p_1,\dots,p_k]$.
\end{proof}

We now prove the main proposition on the $u=1$ specialization. A related multiplication specialization for Negu\c{t}
operators appears in \cite[Lemma~4.6]{DIIP}, with $q=1$ in that
paper's notation.

\begin{proposition}\label{prop:multiplication}
For every symmetric function $f$ whose coefficients are regular at $u=1$,
$\Theta_{a,b}(f)$ is regular there and
\begin{equation}\label{eq:multspecial}
 \Theta_{a,b}(f)\big|_{u=1}
 =\text{multiplication by }\bigl(\Theta_{a,b}(f)1\bigr)\big|_{u=1}.
\end{equation}
In particular, $(\Theta(f)1)_{u=1}$ is multiplicative in $f$.
\end{proposition}
\begin{proof}
Since  $f\in \QQ(q)[[\hbar]][p_1,p_2,\dots]$ by assumption, we have $f\in \QQ(q)[[\hbar]][U_1,U_2,\dots]$ by Lemma~\ref{lem:u-p-transition}.  Because $\Theta_{a,b}(U_n)\in \mathscr{R}$ by Lemma~\ref{lem:qopregular}, $\Theta_{a,b}$ is a $\QQ(q,u)$-algebra homomorphism, and $\mathscr R$ is a $\QQ(q)[[\hbar]]$-algebra, it follows that $\Theta_{a,b}(f)\in \mathscr{R}$.
The proof is thus complete by the second part of Lemma~\ref{lem:rees}. 
\end{proof}

\subsection{The scalar algebra homomorphism}
We combine multiplication specialization with the following scalar
evaluation.  Define the $\QQ(q,u)$-algebra homomorphism
\begin{equation}\label{eq:sign}
 \varepsilon:\Lambda\longrightarrow\QQ(q,u),\qquad
 \varepsilon(p_k)=(-1)^{k-1}.
\end{equation}
In
Section~\ref{subsec:signextraction}, we will use its Schur-coefficient interpretation.
For homogeneous $f\in\QQ(q)[p_1,p_2,\ldots]$ of degree $N$, put
\begin{equation}\label{eq:Phi}
 \Phi(f):=(-1)^{N(b+1)}
       \varepsilon\!\left((\Theta_{a,b}(f)1)|_{u=1}\right),
\end{equation}
and extend linearly across degrees.  The specialization is well defined
by Proposition~\ref{prop:multiplication}.
Here $u=1$ means taking the constant coefficient of the regular
$\hbar$-expansion, not inserting $u=1$ into an uncancelled commutator.

The multiplication specialization allows us to pass from commuting
operators to a multiplicative scalar evaluation.

\begin{proposition}\label{prop:PhiHom}
The map $\Phi:\QQ(q)[p_1,p_2,\ldots]\to\QQ(q)$ is a unital
$\QQ(q)$-algebra homomorphism.
\end{proposition}
\begin{proof}
The map $\Theta$ is an algebra homomorphism into commuting operators.
After specialization, each operator is multiplication by its value on
$1$, by Proposition~\ref{prop:multiplication}.
Consequently $f\mapsto(\Theta(f)1)|_{u=1}$ is an algebra homomorphism
into symmetric functions.  The map $\varepsilon$ is multiplicative,
and the sign $(-1)^{N(b+1)}$ is multiplicative in the degree.
The unit and the coefficient field are fixed.
\end{proof}

\section{The \texorpdfstring{$q$}{q}-difference equation for return paths}\label{sec:Euler}
We use the homomorphism $\Phi$ to prove
$\Pseries(qz;q)=\Aseries(-z;q)\Pseries(z;q)$.
We first derive this equation from two evaluations of $\Phi$, then
prove those evaluations using the compositional shuffle theorem.

\subsection{The two symmetric-function inputs}
With $[z]:=[z^1]$, define the first creation operator
\begin{equation}\label{eq:C1}
 \Cop_1 f[X]:=[z]\left(f[X-(1-q^{-1})/z]H_X(z)\right).
\end{equation}
A composition $\alpha=(\alpha_1,\ldots,\alpha_\ell)\models N$ is an ordered list of positive integers summing to $N$. For $r\ge1$, the general creation operator is
\[
 \Cop_r f[X]:=(-q)^{1-r}[z^r]\left(f[X-(1-q^{-1})/z]H_X(z)\right),
 \qquad \Cop_\alpha:=\Cop_{\alpha_1}\cdots\Cop_{\alpha_\ell}.
\]
In particular, $\Cop_{(1^N)}1=\Cop_1^N 1$ is repeated application of $\Cop_1$ on the seed $1$. We have an identity $e_N=\sum_{\alpha\models N}\Cop_\alpha1$ \cite[(1.6)]{BGLX}. 
The next proposition identifies the path generating functions $P_N(q)$ and $A_N(q)$ as images of these symmetric functions under $\Phi$.

\begin{proposition}\label{prop:seedvalues}
For every $N\ge0$, we have
\begin{equation}\label{eq:Phivalues}
 \Phi(\Cop_1^N1)=P_N(q),\qquad
 \Phi(e_N)=A_N(q).
\end{equation}
\end{proposition}
We defer the proof to Section~\ref{subsec:signextraction}; it uses the compositional shuffle theorem.  We first use the two evaluations
to derive the required functional equation.

\subsection{The repeated first creation operator}
In degree two, \eqref{eq:C1} gives
\[
 \Cop_1^2 1=p_1^2-(1-q^{-1})h_2[X]
 =\frac{1+q^{-1}}2p_1^2-\frac{1-q^{-1}}2p_2.
\]
The general formula below computes $\Cop_1^N 1$.  Its power of $q$ will be what accounts for the normalization of $\Pseries$.

\begin{lemma}\label{lem:seed}
For all $N\ge0$, we have
\begin{equation}\label{eq:seed}
 \Cop_1^N1=q^{-\binom N2}(q)_N h_N[X/(1-q)].
\end{equation}
\end{lemma}
\begin{proof}
For $N=0$, both sides are $1$. To proceed inductively, we first compute $\Cop_1h_{N-1}[X/(1-q)]$. For any alphabet $X=x_1+x_2+\dots$, consider the generating series $H[X]=\sum_{n=0}^\infty h_n[X]=\prod_{i=1}^\infty (1-x_i)^{-1}$, so that $H_X(v)=H[Xv]$, where $v$ is viewed as a monomial alphabet, that is, $p_j[v]:=v^j$. We have $H[X+Y]=H[X]H[Y]$ by \eqref{eq:Hgenerating}.

Let $f[X]=\sum_{N\ge 0} h_N[X/(1-q)]v^N=H[Xv/(1-q)]$, viewed as a formal series in $v$. Applying \eqref{eq:C1} to $f[X]$ coefficientwise in $v$ gives the following computation: 
\[
 \begin{aligned}
 \Cop_1f[X]&=[z]\left(H\left[\frac{(X-(1-q^{-1})/z)v}{1-q}\right]H[Xz]\right)
 =[z]\left(H\left[\frac{Xv}{1-q}+vq^{-1}z^{-1}\right]H[Xz]\right)\\
 &=[z]\left(H\left[\frac{Xv}{1-q}\right]H[vq^{-1}z^{-1}]H[Xz]\right)
 =H\left[\frac{Xv}{1-q}\right]\cdot [z]\frac{H_X(z)}{1-v/(qz)}\\
 &=H\left[\frac{Xv}{1-q}\right]\sum_{j\ge0}q^{-j}h_{j+1}[X]v^j
 =H\left[\frac{Xv}{1-q}\right] \cdot \frac{q}{v}(H[Xv/q]-1)\\
 &=\frac{q}{v}\left(H\left[\frac{Xv}{1-q}+\frac{Xv}{q}\right]-H\left[\frac{Xv}{1-q}\right]\right)\\
 &=\frac{q}{v}\left(H\left[\frac{Xv}{q(1-q)}\right]-H\left[\frac{Xv}{1-q}\right]\right).
 \end{aligned}
\]
The computation above is coefficientwise in $v$, so before coefficient extraction of $z$, it is necessary to work in the iterated Laurent series ring $\QQ(q)[p_1,p_2,\dots]((z))((v))$, with $z$ and $v$ in this order. In particular, $v/z$ is topologically nilpotent, so the only correct expansion of $1/(1-v/(qz))$ is $\sum_{j\ge 0} (v/(qz))^j$.
For $N\ge 1$, comparing coefficients of $v^{N-1}$ gives
\[
 \Cop_1h_{N-1}[X/(1-q)]=q^{1-N}(1-q^N)h_N[X/(1-q)].
\]
Iterating this identity proves \eqref{eq:seed}.
\end{proof}

Applying Proposition~\ref{prop:seedvalues} gives
\begin{equation}\label{eq:PhiH}
 \Phi(h_N[X/(1-q)])=\frac{q^{\binom N2}P_N(q)}{(q)_N}.
\end{equation}
\subsection{From multiplicativity to the scalar equation}
Extend $\Phi$ coefficientwise to formal series in $z$, fixing $z$.
Equation~\eqref{eq:PhiH} and Proposition~\ref{prop:seedvalues} give
\begin{equation}\label{eq:PhiSeries}
 \Pseries(z;q)=\Phi\!\left(\sum_{N\ge0}h_N[X/(1-q)]z^N\right),
 \qquad
 \Aseries(z;q)=\Phi\!\left(\sum_{N\ge0}e_N[X]z^N\right).
\end{equation}
We can therefore apply $\Phi$ to the exponential expressions
\eqref{eq:Hgenerating}--\eqref{eq:Egenerating}.
Writing $r_k=\Phi(p_k)$ and using multiplicativity yields

\begin{align}
 \Pseries(z;q)&=\exp\left(\sum_{k\ge1}\frac{r_k z^k}{k(1-q^k)}\right),\label{eq:Hexponential}\\
 \Aseries(z;q)&=\exp\left(\sum_{k\ge1}\frac{(-1)^{k-1}r_k z^k}{k}\right).
 \label{eq:Aexponential}
\end{align}
Both expressions are formal in $z$ over $\QQ(q)$, so every coefficient uses only finitely many terms.

Comparing these two exponentials at $z$ and at $qz$ produces the promised $q$-difference equation for the return-path series.
\begin{theorem}\label{thm:Euler}
The return-path series \eqref{eq:pathseries} satisfies
\begin{equation}\label{eq:HEuler}
 \Pseries(qz;q)=\Aseries(-z;q)\Pseries(z;q),\qquad\Pseries(0;q)=1.
\end{equation}
\end{theorem}
\begin{proof}
Subtract the logarithm of \eqref{eq:Hexponential} from its value at $qz$:
\[
 \sum_{k\ge1}\frac{r_k(q^k-1)z^k}{k(1-q^k)}
 =-\sum_{k\ge1}\frac{r_kz^k}{k}
 =\log\Aseries(-z;q).
\]
Exponentiation proves the result.
\end{proof}

By \eqref{eq:HP}, the same equation holds for $\Hseries$.
Therefore, its coefficients are the $G_N$ of \eqref{eq:solutionrecurrence}.
Substituting $G_N=q^{\gamma_N}F_N/(q)_N$ in that recurrence,
multiplying by $(q)_N$, and using $(q)_N/(1-q^N)=(q)_{N-1}$ gives
\begin{equation}\label{eq:finiteEuler}
 q^{\gamma_N}F_N(q)=\sum_{k=1}^N(-1)^{k+1}A_k(q)
       q^{\gamma_{N-k}}F_{N-k}(q)\frac{(q)_{N-1}}{(q)_{N-k}}
 \qquad(N\ge1).
\end{equation}
The quotient is a finite product, so every summand is a polynomial.

\begin{remark}
For $N=1$, \eqref{eq:finiteEuler} says $F_1=A_1$, the area/dinv symmetry of the rational $q,t$-Catalan numbers. For $N>1$, the right side contains terms of degree below $\gamma_N$, but their sum is divisible by $q^{\gamma_N}$. The symmetric-function proof establishes this divisibility.
 Section~\ref{subsec:runningrecurrence} illustrates the cancellation at $N=2$.
\end{remark}

\subsection{Parking functions and the shuffle input}\label{subsec:parking}
We turn to the proof of Proposition~\ref{prop:seedvalues}.
The required external input is the compositional rational shuffle
identity formulated in \cite[(3.7)]{BGLX} and proved in
\cite[Section~6]{Mellit}.
Fix $N\ge1$, and put $m=aN$ and $n=bN$.
  If the $i$th north step starts at $(x_i,i-1)$, its rank is
\begin{equation}\label{eq:PFrank}
 R_i:=m(i-1)-nx_i+\frac{x_i}{m+1}.
\end{equation}
The ranks are distinct.  The perturbation separates equal integral ranks
without reversing unequal ones, and ranks increase up each column.
Order the north steps by their ranks
$\rho_1<\cdots<\rho_n$.  The rank window of $P$ is
\begin{equation}\label{eq:rankwindow}
 \mathcal W(P):=\{(i,j):1\le i<j\le n,\quad \rho_j-\rho_i<m\},
 \qquad \max\tdinv(P):=|\mathcal W(P)|.
\end{equation}

A {\it parking function} $\pi$ on a path $P\in\mathcal D_{m,n}$ is a labelling of north steps by $1,\ldots,n$ that is increasing from bottom to top in each column. For a parking function $\pi$ on $P$, let $w_i$ be the label of the step
of rank $\rho_i$.  Its decreasing-rank reading word is
$\sigma_\pi=(w_n,\ldots,w_1)$, and its inverse descent set is
\[
 \ides(\pi):=\{r\in\{1,\ldots,n-1\}:r+1\text{ precedes }r
                   \text{ in }\sigma_\pi\}.
\]
Define its temporary dinv by
\begin{equation}\label{eq:tdinvwindow}
 \tdinv(\pi):=\#\{(i,j)\in\mathcal W(P):w_i<w_j\}.
\end{equation}
Therefore, $\tdinv(\pi)\le\max\tdinv(P)$ directly from the definition.
Equivalently, if $R(r)$ is the rank of the step labelled $r$, this counts
$r<s$ with $R(r)<R(s)<R(r)+m$.

On each $P$, there is a distinguished \emph{rank-labelled parking function} $\pi_P^{\mathrm{rk}}$ that assigns
label $i$ to the step of rank $\rho_i$.  It is legal because ranks
increase up each column.  Its reading word is $(n,n-1,\ldots,1)$, and
\begin{equation}\label{eq:ranklabelled}
 \ides(\pi_P^{\mathrm{rk}})=\{1,\ldots,n-1\},\qquad
 \tdinv(\pi_P^{\mathrm{rk}})=|\mathcal W(P)|.
\end{equation}
This proves that \eqref{eq:rankwindow} is the maximum temporary dinv
among all legal labellings, justifying the notation.  Moreover,
$\pi_P^{\mathrm{rk}}$ is the unique labelling with full inverse descent set,
since that set forces the reading word to be $(n,n-1,\ldots,1)$.

For a parking function $\pi$ with underlying path $P_\pi$, set
\begin{equation}\label{eq:PFdinv}
 \dinv(\pi):=\dinv(P_\pi)+\tdinv(\pi)-\max\tdinv(P_\pi).
\end{equation}
In particular, $\dinv(\pi_P^{\mathrm{rk}})=\dinv(P)$.
Here $\dinv(P)$ is the arm-leg count \eqref{eq:hook}.
These are the statistics of \cite[Algorithm~6.1]{BGLX}, with $u$ in place
of its area parameter $t$.  See also \cite{Haglund} for background.

\begin{example}[A parking function with partial inverse descent set]
Take $(a,b,N)=(2,3,2)$ and the path
$P=\mathtt{NNENENNENE}$ in the $4\times6$ rectangle.
Label its north steps in chronological order by $1,4,2,3,6,5$.
The labels increase within each column.
Figure~\ref{fig:parking} and Table~\ref{tab:parking} show the path and
its perturbed ranks.
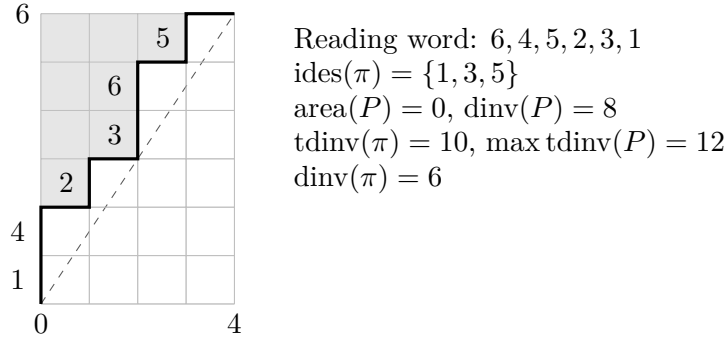
\begin{figure}[!htbp]
\centering
\begin{tikzpicture}[x=.64cm,y=.64cm,font=\small]
 \draw[step=1,black!20,very thin] (0,0) grid (4,6);
 \draw[dashed,black!65] (0,0)--(4,6);
 \fill[black!10] (0,2) rectangle (1,6);
 \fill[black!10] (1,3) rectangle (2,6);
 \fill[black!10] (2,5) rectangle (3,6);
 \draw[step=1,black!25,very thin] (0,0) grid (4,6);
 \draw[very thick] (0,0)--(0,2)--(1,2)--(1,3)--(2,3)
          --(2,5)--(3,5)--(3,6)--(4,6);
 \foreach \x/\y/\l in {0/0/1,0/1/4,1/2/2,2/3/3,2/4/6,3/5/5}
  {\node[anchor=east] at (\x-.10,\y+.5) {$\l$};}
 \node[below] at (0,0) {$0$};
 \node[below] at (4,0) {$4$};\node[left] at (0,6) {$6$};
 \node[anchor=west,align=left] at (5,4)
  {Reading word: $6,4,5,2,3,1$\\
   $\ides(\pi)=\{1,3,5\}$\\
   $\area(P)=0$, $\dinv(P)=8$\\
   $\tdinv(\pi)=10$, $\max\tdinv(P)=12$\\
   $\dinv(\pi)=6$};
\end{tikzpicture}
\caption{An above-diagonal parking function.  The north-step labels
are written to the left of the corresponding steps.
The shaded Ferrers diagram has eight boxes, all counted by $\dinv(P)$.}
\label{fig:parking}
\end{figure}
\par
\begin{table}[!htbp]
\centering\renewcommand{\arraystretch}{1.18}
\begin{tabular}{|c|c|c|c|c|}
\hline
North step $i$ & Foot $(x_i,i-1)$ & Label & Rank $R_i$ & Rank index\\\hline
$1$&$(0,0)$&$1$&$0$&$1$\\\hline
$2$&$(0,1)$&$4$&$4$&$5$\\\hline
$3$&$(1,2)$&$2$&$2+\tfrac15$&$3$\\\hline
$4$&$(2,3)$&$3$&$\tfrac25$&$2$\\\hline
$5$&$(2,4)$&$6$&$4+\tfrac25$&$6$\\\hline
$6$&$(3,5)$&$5$&$2+\tfrac35$&$4$\\\hline
\end{tabular}
\caption{Labels and perturbed ranks in Figure~\ref{fig:parking}.  The rank index is the position in ascending rank order.}
\label{tab:parking}
\end{table}
\par
The ascending ranks and their labels are
\[
 (\rho_1,\ldots,\rho_6)=
 \left(0,\tfrac25,2+\tfrac15,2+\tfrac35,4,4+\tfrac25\right),
 \qquad (w_1,\ldots,w_6)=(1,3,2,5,4,6).
\]
Therefore, the reading word is $(6,4,5,2,3,1)$ and
$\ides(\pi)=\{1,3,5\}$.
There are twelve pairs in $\mathcal W(P)$; for example, $(1,5)\notin \mathcal W(P)$, since $\rho_5-\rho_1=4-0=m$, violating the strict upper endpoint in \eqref{eq:rankwindow}.  All but two, $(2,3)$ and $(4,5)$, have their labels in increasing order.
Consequently $\tdinv(\pi)=12-2=10$ and $\dinv(\pi)=8+10-12=6$.
The rank-labelled parking function instead has chronological labels
$(1,5,3,2,6,4)$ and its tdinv counts all twelve pairs.
\end{example}

For $S\subseteq\{1,\ldots,n-1\}$, the fundamental quasisymmetric function is
\[
 F_{n,S}[X]:=\sum_{\substack{i_1\le\cdots\le i_n\\j\in S\Rightarrow i_j<i_{j+1}}}
 x_{i_1}\cdots x_{i_n}.
\]
Individual fundamental functions need not be symmetric.
For a composition $\alpha=(\alpha_1,\ldots,\alpha_\ell)\models N$,
let $\operatorname{PF}_{aN,bN}^{\alpha}$ be the set of parking functions
whose underlying paths meet the diagonal exactly at
\[
 (ak_j,bk_j),\qquad k_j=\alpha_1+\cdots+\alpha_j,
 \quad 0\le j\le\ell,
\]
where $k_0=0$.  Therefore, $\alpha$ records successive return distances in
units of $(a,b)$, in the order the path is traversed.
For example, a path in $\operatorname{PF}_{5a,5b}^{(3,1,1)}$
returns at $(3a,3b)$ and $(4a,4b)$ but not at $(a,b)$ or $(2a,2b)$.

The following theorem, conjectured by Bergeron--Garsia--Leven--Xin \cite[(3.7)]{BGLX} and proved by Mellit \cite[Section~6]{Mellit}, evaluates the slope operators on compositional seeds as weighted sums over parking functions with prescribed diagonal returns.
\begin{theorem}[Rational compositional shuffle theorem]\label{thm:shuffleinput}
With the preceding conventions, for $N\ge1$ and every composition $\alpha\models N$, we have
\begin{equation}\label{eq:shuffleinput}
 (-1)^{N(b+1)}\Theta_{a,b}(\Cop_\alpha1)1
 =\sum_{\pi\in\operatorname{PF}^{\alpha}_{aN,bN}}
 q^{\dinv(\pi)}u^{\area(\pi)}F_{bN,\ides(\pi)}[X].
\end{equation}
\end{theorem}

\subsection{From the shuffle identity to the two evaluations}
\label{subsec:signextraction}
To prove Proposition~\ref{prop:seedvalues}, we extract the column Schur coefficient from the shuffle identity.  For a homogeneous symmetric function $f$ of degree $n$, note that
\begin{equation}\label{eq:epsilonSchur}
 \varepsilon(f)=\langle s_{(1^n)},f\rangle,
\end{equation}
where $\langle -,-\rangle$ is the Hall inner product under which the Schur basis is orthonormal.
Indeed, for the algebra involution $\omega(p_k)=(-1)^{k-1}p_k$, we have
$\varepsilon(f)=\omega(f)[1]$.
The identities $\omega(s_\lambda)=s_{\lambda'}$ and
$s_\mu[1]=0$ unless $\mu$ has one row prove
\eqref{eq:epsilonSchur} \cite[Chapter~I, Sections~3--4]{Macdonald}.

To read this coefficient in a fundamental expansion, use Gessel's
identity \cite{Gessel}
\[
 s_\lambda=\sum_{T\in\operatorname{SYT}(\lambda)}F_{n,\operatorname{Des}(T)}.
\]
Here $\operatorname{SYT}(\lambda)$ is the set of standard Young tableaux
of shape $\lambda$, whose entries increase along rows and columns.
The descent set $\operatorname{Des}(T)$ consists of the $j$ for which $j+1$
lies in a lower row than $j$.
A standard tableau has every possible descent only when it is the
unique column tableau.  Consequently, for a symmetric sum of
fundamental functions, its $F_{n,\{1,\ldots,n-1\}}$ coefficient is
exactly its $s_{(1^n)}$ coefficient.

For each underlying path $P$, full inverse descent selects the
rank-labelled parking function $\pi_P^{\mathrm{rk}}$ of
\eqref{eq:ranklabelled}, with contribution $q^{\dinv(P)}u^{\area(P)}$.
Applying this extraction to the all-ones composition and to the sum over
compositions gives the following unlabelled evaluations.

\begin{proposition}\label{prop:unlabelledshuffle}
For $N\ge1$, we have
\begin{align}
 (-1)^{N(b+1)}\varepsilon(\Theta(\Cop_1^N1)1)
   &=\sum_{P\in\Rpaths}q^{\dinv(P)}u^{\area(P)},
       \label{eq:shuffle-return}\\
 (-1)^{N(b+1)}\varepsilon(\Theta(e_N)1)
   &=\sum_{P\in\mathcal D_{aN,bN}}q^{\dinv(P)}u^{\area(P)}.
       \label{eq:shuffle-all}
\end{align}
\end{proposition}
\begin{proof}
Apply \eqref{eq:epsilonSchur} and the full-descent argument to
Theorem~\ref{thm:shuffleinput}.
The composition $(1^N)$ forces a return after each coprime block.
Summing over all compositions, and using
$e_N=\sum_{\alpha\models N}\Cop_\alpha1$, counts all paths.

\end{proof}

\begin{proof}[Proof of Proposition~\ref{prop:seedvalues}]
For $N=0$, both evaluations are $\Phi(1)=1$.
For $N\ge1$, specialize \eqref{eq:shuffle-return} at $u=1$ to obtain
\[
 \Phi(\Cop_1^N1)=\sum_{P\in\Rpaths}q^{\dinv(P)}=P_N(q).
\]
For the second evaluation, $\Theta(e_N)1$ is invariant under
exchanging $q,u$: the basic operators and their commutator recursion
depend symmetrically on these parameters, and the generators $U_k$
depend only on $qu$.  The seed $e_N$ is parameter independent.
The shuffle expansion makes both specializations legitimate.
Its sign coefficient at $(q,1)$ therefore equals its value at $(1,q)$,
which is $A_N(q)$ by \eqref{eq:shuffle-all}.
\end{proof}

\begin{remark}
The exchange argument shows $A_N(q)=\sum_{P\in\mathcal D_{aN,bN}}q^{\dinv(P)}$. If one adopted this as the definition of $A_N(q)$ instead of the one using areas, then the proof of Proposition~\ref{prop:seedvalues} would not require the exchange argument. However, the area/dinv symmetry in $\mathcal D_{aN,bN}$ is unavoidable, because the proof of Theorem~\ref{thm:cofactor} requires the area definition of $A_N(q)$ via Lemma~\ref{lem:walkarea}.
\end{remark}

\subsection{A rank-two specialization}\label{subsec:runningrecurrence}
The case $(a,b)=(2,3)$ illustrates the cancellations in the finite
recurrence.  Its unrestricted area polynomial at rank two is
\begin{equation}\label{eq:exampleA2}
 A_2(q)=1+3q+4q^2+4q^3+4q^4+3q^5+2q^6+q^7+q^8.
\end{equation}
It follows by enumerating
$y_0\le y_1\le y_2\le y_3=6$ with
$y_x\ge\lceil3(x+1)/2\rceil$ and weight $q^{y_0+y_1+y_2-10}$.
Since $A_1=F_1=1+q$, \eqref{eq:finiteEuler} becomes
\begin{equation}\label{eq:exampleEuler}
 q^5F_2=(1+q)^2-(1-q)A_2=q^5+q^6+q^7+q^9.
\end{equation}
This agrees with \eqref{eq:exampleF2}.
The unrestricted family has $A_2(1)=23$ paths, whereas the return
family has $F_2(1)=4$.  The two terms in the middle expression of \eqref{eq:exampleEuler} each contain powers below $q^5$, which cancel in their difference.

\section{The rank graph and disjoint cycle series}\label{sec:graph}
We next encode rational Dyck paths as walks in a weighted graph and
define a generating series for disjoint cycles.  The cycle lengths and
weights established here will be used in both the determinant argument
and the cylinder bijection.

Let $a,b>1$ be relatively prime and $d=a+b$. For indeterminate variables $q,s$, let $\Ggraph$ be a weighted directed graph with vertices $r\in\NN$ and directed edges
\begin{equation}\label{eq:rankgraphedges}
 \begin{aligned}
 r&\longrightarrow r+a &&(r\ge0),
     &\text{weight }&s,\\
 r&\longrightarrow r-b &&(r\ge b),
     &\text{weight }&s q^{\lfloor(r-b)/a\rfloor}.
 \end{aligned}
\end{equation}
These are the rank increments of a north step and an east step for the rank function
$ay-bx$, respectively. 
A walk on the graph is allowed to repeat vertices; its weight is the product of all edge weights. A cycle is a directed simple cycle
with distinct vertices apart from the coinciding endpoints.
A collection of cycles is vertex-disjoint when no vertex belongs to two cycles.

A walk from $0$ to $0$ of length $dN$ has $bN$ north edges and $aN$
east edges.  Reading these in their given order yields an above-diagonal
path in $\mathcal D_{aN,bN}$, because every visited rank is nonnegative.
Conversely, the vertex ranks of such a path give a unique graph walk.

The edge weights recover the area of the corresponding path.

\begin{lemma}\label{lem:walkarea}
The weight of a length-$dN$ walk from $0$ to $0$ is
$s^{dN}q^{\area(P)}$, where $P$ is its associated above-diagonal path.
\end{lemma}
\begin{proof}
An east step of $P$ ends at $(x+1,y_x)$, with rank
$r=ay_x-b(x+1)$.  Its edge contributes
\[
 \lfloor r/a\rfloor=y_x-\left\lceil\frac{b(x+1)}a\right\rceil
\]
to the exponent of $q$.  Summing over $0\le x<aN$ gives
\eqref{eq:introarea}, and the $dN$ edges contribute $s^{dN}$.
\end{proof}

The residue classes modulo $d$ determine the structure of every simple
cycle.  They also determine how disjoint cycles occupy the graph.

\begin{lemma}\label{lem:cyclelength}
Every cycle has $d$ vertices, $b$ edges of increment $+a$, and $a$
edges of increment $-b$.
A vertex-disjoint collection is determined by its occupied vertices.
A collection of $N$ cycles occupies $N$ vertices in each residue
class modulo $d$.
\end{lemma}
\begin{proof}
Let $\sigma$ send each occupied vertex to its successor.
Both steps increase its residue modulo $d$ by $a$.
For occupied $x<y$ of the same residue, $y-x\ge d$, and
$\sigma(y)-\sigma(x)$ is one of $y-x-d,y-x,y-x+d$.
It is nonnegative and cannot vanish because $\sigma$ is injective.
Thus $\sigma$ maps each residue class order-preservingly onto the next.
Since $\gcd(a,d)=1$, these residue classes form one orbit and have
equal sizes.  Their $j$th lowest vertices form one cycle of length
$d$.  This also determines every edge from the occupied set.
If a cycle has $u$ north edges and $v$ east edges, then
$au=bv$ and $u+v=d$, so $u=b$ and $v=a$.
\end{proof}

Write $w(C)$ for a cycle's $q$-exponent, so the weight of a (length-$d$) cycle is $s^d q^{w(C)}$.
For a finite cycle collection $\mathcal F$, we also write
\[
 S(\mathcal F):=\sum_{C\in\mathcal F}S(C),\qquad
 w(\mathcal F):=\sum_{C\in\mathcal F}w(C).
\]

The next formula determines both its minimum weight and the effect of
translating its vertices.
\begin{lemma}\label{lem:cycleweight}
If $S(C)$ is the sum of the vertices of $C$, then
\begin{equation}\label{eq:cycleweight}
 w(C)=\frac{S(C)}d-\frac{d-1}2.
\end{equation}
Thus, translating $C$ upward by one increases $w(C)$ by one.
For a collection $\mathcal{F}$ of $N$ disjoint cycles occupying vertices
$r+d k_{r,j}$, where $0\le r<d$ and
$0\le k_{r,1}<\cdots<k_{r,N}$,
\begin{equation}\label{eq:totalweight}
 w(\mathcal F)=\sum_{r,j}k_{r,j}\ge d\binom N2=\gamma_N.
\end{equation}
\end{lemma}
\begin{proof}
The weight is the sum of the floors of the east-edge ending ranks.
We first compute the sum of those ranks and then subtract their
fractional parts.

Let $U$ be the sum of the starting ranks of north edges $u\to u+a$, and let $E$
be the sum of the ending ranks of east edges $r+b\to r$.  There are $b$ north
edges and $a$ east edges.  Consider
\[
 0=\sum_{\text{north starts }u}\bigl((u+a)^2-u^2\bigr)
   +\sum_{\text{east ends }r}\bigl(r^2-(r+b)^2\bigr)
   =2aU+ba^2-2bE-ab^2;
\]
the middle expression vanishes because the squared rank of each vertex is cancelled out.
The starting ranks of the east edges sum to $E+ab$, so
$S(C)=U+E+ab$.  Eliminating $U$ from these two identities gives
\[
 \frac Ea=\frac{S(C)}d-\frac b2.
\]
Between consecutive east edges, north edges do not change the rank
modulo $a$.  Each east edge subtracts $b$, and $\gcd(a,b)=1$.
The $a$ east-edge ending ranks therefore have residues
$0,1,\ldots,a-1$ modulo $a$, in some order.  Hence
\[
 w(C)=\sum_{\text{east ends }r}\lfloor r/a\rfloor
      =\frac Ea-\frac{0+1+\cdots+(a-1)}a
      =\frac{S(C)}d-\frac{d-1}2.
\]
This proves \eqref{eq:cycleweight}.  

Translating all $d$ vertices
upward by one adds $d$ to $S(C)$ and thus adds one to $w(C)$.
For a collection of $N$ cycles, the sum of vertices is
\[
 S(\mathcal F)=\sum_{r=0}^{d-1}\sum_{j=1}^N(r+dk_{r,j})
 =N\frac{d(d-1)}2+d\sum_{r,j}k_{r,j},
\]
so
\[
 w(\mathcal F)=\frac{S(\mathcal F)}d-N\frac{d-1}2
 =\sum_{r,j}k_{r,j}.
\]
Finally, the strict inequalities $0\le k_{r,1}<\dots<k_{r,N}$ force $\sum_j k_{r,j}\ge 0+1+\dots+(N-1)=\binom{N}{2}$, so $\sum_{r,j}k_{r,j}\ge \gamma_N.$
\end{proof}

We can now define the disjoint cycle series directly:
\[ \mathcal{DC}_N:=\{\text{vertex-disjoint collections of $N$ cycles in $\Ggraph$}\}, \quad \mathcal{DC}:=\bigsqcup_{N\ge 0} \mathcal{DC}_N.\]
\begin{equation}\label{eq:cycleseries}
\begin{gathered}
D_N(q):=\sum_{\mathcal F\in \mathcal{DC}_N} q^{w(\mathcal F)}, \qquad  w(\mathcal F):=\sum_{C\in\mathcal F}w(C),\\
 \Dseries(z;q):=\sum_{N\ge 0} D_N(q)z^N=\sum_{\mathcal F\in \mathcal{DC}}z^{|\mathcal F|}q^{w(\mathcal F)},
 \end{gathered}
\end{equation}
where $\abs{\mathcal F}$ denotes the number of cycles rather than the number of vertices occupied.

The cycle-weight formula bounds the vertices that can contribute to any fixed coefficient.  This gives a well-defined generating series.

\begin{proposition}\label{prop:detlimit}
The series $D_N(q)$ is well defined and divisible by $q^{\gamma_N}$. The series $\Dseries(z;q)$ is well defined in $\NN[[q]][[z]]$ with constant term $1$.
\end{proposition}
\begin{proof}
From \eqref{eq:cycleweight}, there are only finitely many cycles $C$ with bounded $w(C)$, so there are only finitely many contributions to a given coefficient of $D_N(q)$. The divisibility by $q^{\gamma_N}$ follows from \eqref{eq:totalweight}.
\end{proof}

\section{The cycle-series equation from finite determinants}\label{sec:determinant}
We prove the $q$-difference equation for $\Dseries$ using determinantal interpretations
of finite graph truncations.  
For $H\ge0$, let
\begin{equation}\label{eq:detcycles}
 \Dseries_H(z;q):=\sum_{\mathcal F\subseteq[0,H]}
                z^{|\mathcal F|}q^{w(\mathcal F)},
 \qquad
 \Aseries_H(z;q):=\sum_{N\ge0}A_{N,H}(q)z^N,
\end{equation}
where the first sum restricts the occupied vertices of each collection
to $[0,H]$, and $A_{N,H}$ counts length-$dN$ excursions (i.e., walks from $0$ to $0$)
inside this interval, with the $q$-weight of Section~\ref{sec:graph}.
The first series is a polynomial in $q,z$ with constant term $1$.
Proposition~\ref{prop:detlimit} gives
$\Dseries_H\to\Dseries$ coefficientwise in $q,z$.
Also $A_{N,H}=A_N$ for $H\ge abN$, since a length-$dN$ excursion
cannot reach a rank above $abN$. Hence, $\Aseries_H\to \Aseries$ coefficientwise in $q,z$.

Let $L_H(s,q)$ be the weighted adjacency matrix of this finite graph,
with $(i,j)$-entry equal to the weight of $i\to j$; if there is no edge, the corresponding entry is zero.
The following lemma says the permutation expansion of the determinant counts vertex-disjoint
cycles, with one minus sign for each cycle.

\begin{lemma}\label{lem:detpositive} Assuming the notation above, we have
\begin{equation}\label{eq:detpositive}
 \det(I-L_H(s,q))=\Dseries_H(-s^d;q).
\end{equation}
\end{lemma}
\begin{proof}
A nontrivial permutation cycle contributes only if its edges lie in the
graph.  Thus determinant terms correspond to disjoint simple cycles,
as in the transfer-matrix method \cite[Section~4.7]{StanleyEC1}.
A cycle of length $\ell$ contributes permutation sign
$(-1)^{\ell-1}$ and entry sign $(-1)^\ell$, whose product is $-1$.
By Lemma~\ref{lem:cyclelength}, $\ell=d$.
A collection of $N$ cycles therefore contributes
$(-s^d)^Nq^{w(\mathcal F)}$.
\end{proof}

Deleting vertex $0$ and subtracting one from the remaining vertex labels
produces the graph with cutoff $H-1$.  The cycle-weight formula determines
the corresponding change in its determinant.

\begin{theorem}\label{thm:cofactor}
For $H\ge1$, we have
\begin{equation}\label{eq:finitecofactor}
 \Aseries_H(-z;q)\Dseries_H(z;q)=\Dseries_{H-1}(qz;q).
\end{equation}
Consequently, we have
\begin{equation}\label{eq:Dfunctional}
 \Dseries(qz;q)=\Aseries(-z;q)\Dseries(z;q).
\end{equation}
\end{theorem}
\begin{proof}
Since every matrix entry of $L_H$ is divisible by $s$, the inverse
$(I-L_H)^{-1}=\sum_{m\ge0}L_H^m$ is entrywise a formal series in $\ZZ[[q,s]]$.
Its $(0,0)$ entry counts excursions and equals
$\Aseries_H(s^d;q)$ by Lemma~\ref{lem:walkarea}.
Cramer's rule expresses this entry as the principal cofactor deleting
vertex $0$, divided by $\det(I-L_H)$. Note that $\det(I-L_H)$ is invertible in $\ZZ[[q,s]]$ because it has constant term $1$.

The cofactor matrix is $I-L_{[1,H]}(s,q)$, with the weighted adjacency matrix $L_{[1,H]}(s,q)$ of $\Ggraph$ restricted to the vertex set $[1,H]$.
Translate its vertices down by one.
By Lemma~\ref{lem:cycleweight}, a collection of $N$ cycles has its
original monomial weight equal to $q^N$ times the translated monomial weight. By proof of Lemma~\ref{lem:detpositive}, we have
\[ \det(I-L_{[1,H]}(s,q))=\sum_{\mathcal F\subseteq [1,H]} (-s^d)^{\abs{\mathcal F}}q^{w(\mathcal{F})}=\sum_{\mathcal F\subseteq [0,H-1]} (-s^d)^{\abs{\mathcal F}}q^{\abs{\mathcal F}}q^{w(\mathcal{F})}=\Dseries_{H-1}(-qs^d;q).\]

Consequently, the Cramer's rule for the $(0,0)$ entry reads:
\[ \Aseries_H(s^d;q)=\frac{\det(I-L_{[1,H]}(s,q))}{\det(I-L_H(s,q))}=\frac{\Dseries_{H-1}(-qs^d;q)}{\Dseries_H(-s^d;q)}.\]
This is an identity on the subring $\ZZ[[q,s^d]]$.
Replacing $s^d$ by $-z$ proves \eqref{eq:finitecofactor}.

Since every factor of \eqref{eq:finitecofactor} has a coefficientwise limit in $\ZZ[[q,z]]$, passing to the limit proves \eqref{eq:Dfunctional}.
\end{proof}
Thus $\Dseries$ satisfies the same equation as $\Pseries$ and $\Hseries$.
It remains to identify its cycle coefficients with the bounded-cylinder
coefficients of $\Bseries$.

\section{An explicit disjoint-cycle/cylinder bijection}\label{sec:bijection}
We identify $\Dseries$ with $\Bseries$ by converting occupied vertices
into partitions and placing their conjugates along the diagonals of a
cylinder.  The number of cycles becomes the largest-entry bound, and
the cycle weight becomes volume after subtracting $\gamma_N$. Recall the set $\mathcal{DC}_N$ from \eqref{eq:cycleseries}.

\begin{proposition}
\label{prop:cyclecylinder}
For each $N\ge0$, there is a bijection between disjoint cycle collections $\mathcal F\in \mathcal{DC}_N$ and cylindric partitions
$\boldsymbol\lambda$ of profile \eqref{eq:balanced} with largest entry
at most $N$.  It satisfies
\[
 w(\mathcal F)=\gamma_N+|\boldsymbol\lambda|.
\]
Consequently, we have
\begin{equation}\label{eq:cyclecylindertarget}
 D_N(q)=\sum_{\mathcal F\in \mathcal{DC}_N}q^{w(\mathcal F)}
        =q^{\gamma_N}C_{\cc,\le N}(q),
 \qquad \Dseries=\Bseries.
\end{equation}
\end{proposition}
\subsection{Interlacing partitions from residue runners}\label{subsec:runnermap}
We first describe an intermediate bijection with interlacing partitions.
\begin{definition}
    For partitions $\alpha,\beta$, write $\alpha\prec \beta$ when the skew partition $\beta/\alpha$ is a skew row, that is, $\beta/\alpha$ has at most one box in each column. Equivalently, we have
    \begin{equation}\label{eq:interlacing-inequality}
        \beta_j\ge \alpha_j\ge \beta_{j+1}\text{ for }j\ge 1.
    \end{equation}

    Write $\nu_{t,j}$ for the $j$-th part of $\nu_t$. A sequence of partitions $(\nu_0,\dots,\nu_{d-1})$ is \emph{interlacing} with \emph{interlacing code} $(\eta_0,\dots,\eta_{d-1})\in \{0,1\}^d$ if 
    \begin{equation}\label{eq:interlacing}\nu_t\prec \nu_{t+1} \text{ if }\eta_t=1, \qquad \nu_t\succ \nu_{t+1} \text{ if }\eta_t=0, \qquad 0\leq t<d.\end{equation}
    The index $t$ is understood cyclically modulo $d$.
\end{definition}

Consider the balanced interlacing code
\begin{equation}
    \label{eq:balancedinterlacing}
    \eta_t:=\left\lfloor\frac{a(t+1)}d\right\rfloor-
          \left\lfloor\frac{at}d\right\rfloor,\qquad 0\le t<d.
\end{equation}
For example, for $(a,b)=(3,5)$, $\eta=(0,0,1,0,0,1,0,1)$. The interlacing code has $a$ ones and $b$ zeros.
The residue-runner construction gives the following intermediate form of the cycle--cylinder bijection.

\begin{proposition}\label{prop:cycleinterlacing}
    For $N\ge 0$, there is a bijection between $\mathcal F\in \mathcal{DC}_N$ and interlacing sequences of partitions $(\nu_0,\dots,\nu_{d-1})$ with interlacing code \eqref{eq:balancedinterlacing} and parts bounded by $N$. Moreover, we have
    \begin{equation}\label{eq:cyclevolume}
     w(\mathcal F)=\gamma_N+\sum_{t=0}^{d-1}|\nu_t|.
    \end{equation}
\end{proposition}

To construct the bijection, let $\mathcal F\in \mathcal{DC}_N$. Arrange the occupied vertices into columns (runners) according to their residue classes modulo $d$. More precisely, for a column index $0\leq t<d$, define the corresponding residue class by
\begin{equation}\label{eq:runners}
 r_t:=at\bmod d\in\{0,\ldots,d-1\},\qquad 0\le t<d.
\end{equation}
For a row index $k\ge 0$, let the rank of the $t$-th column and $k$-th row be $r_t+dk$. Hence each nonnegative rank is assigned a unique row and column index. 
Since $+a\equiv-b\equiv a\pmod d$, each step moves from one column to the next, with column indices understood cyclically modulo $d$.

\begin{figure}[!htbp]
\centering
\begin{tikzpicture}[x=.44cm,y=.44cm,font=\scriptsize,>=Stealth]
\foreach \t/\x/\res in {0/0/0,1/3/2,2/6/4,3/9/1,4/12/3,5/15/0}{
  \node at (\x+.6,11.65) {$t=\t$};
  \foreach \k in {0,1,2,3,4}{
    \pgfmathtruncatemacro{\rank}{\res+5*\k}
    \fill[black!3] (\x,5.0+1.28*\k) rectangle ++(1.2,1.15);
    \draw[black!28] (\x,5.0+1.28*\k) rectangle ++(1.2,1.15);
    \node[black!55] at (\x+.6,5.575+1.28*\k) {$\rank$};
  }
}
\draw[dashed,black!65,thick] (14.5,-0.48) rectangle (16.7,12.17);
\foreach \k in {0,1,2,3,4}{\node[anchor=east] at (-.35,5.575+1.28*\k) {$k=\k$};}
\foreach \x/\k/\v in {0/4/20,0/2/10,0/1/5,
                         3/3/17,3/2/12,3/0/2,
                         6/2/14,6/1/9,6/0/4,
                         9/3/16,9/1/6,9/0/1,
                         12/3/18,12/1/8,12/0/3,
                         15/4/20,15/2/10,15/1/5}{
  \fill[black!18] (\x,5.0+1.28*\k) rectangle ++(1.2,1.15);
  \draw[black!78] (\x,5.0+1.28*\k) rectangle ++(1.2,1.15);
  \node at (\x+.6,5.575+1.28*\k) {$\v$};
}
\foreach \xa/\ka/\xb/\kb in {0/4/3/3,3/3/6/2,6/2/9/3,9/3/12/3,12/3/15/4,
                               0/2/3/2,3/2/6/1,6/1/9/1,9/1/12/1,12/1/15/2,
                               0/1/3/0,3/0/6/0,6/0/9/0,9/0/12/0,12/0/15/1}{
  \draw[->,black!78] (\xa+1.2,5.575+1.28*\ka)--(\xb,5.575+1.28*\kb);
}
\node[anchor=east] at (-.35,4.15) {$r_t$};
\foreach \x/\v in {0/0,3/2,6/4,9/1,12/3,15/0}{\node at (\x+.6,4.15) {$\v$};}
\node[anchor=east] at (-.35,3.25) {$k_t$};
\foreach \x/\v in {0/{(4,2,1)},3/{(3,2,0)},6/{(2,1,0)},9/{(3,1,0)},12/{(3,1,0)},15/{(4,2,1)}}{\node at (\x+.6,3.25) {$\v$};}
\node[anchor=east] at (-.35,1.85) {$\mu_t$};
\foreach \x in {0,15}{
  \draw (\x+.25,2.35) rectangle ++(.34,-.34); \draw (\x+.59,2.35) rectangle ++(.34,-.34);
  \draw (\x+.25,2.01) rectangle ++(.34,-.34); \draw (\x+.25,1.67) rectangle ++(.34,-.34);}
\foreach \r in {0,1}{\draw (3.43,2.35-.34*\r) rectangle ++(.34,-.34);}
\node at (6.6,1.95) {$\varnothing$};
\foreach \x in {9,12}{\draw (\x+.43,2.18) rectangle ++(.34,-.34);}
\node[anchor=east] at (-.35,.15) {$\nu_t$};
\foreach \x in {0,15}{
  \foreach \c in {0,1,2}{\draw (\x+.08+.34*\c,.62) rectangle ++(.34,-.34);}
  \draw (\x+.08,.28) rectangle ++(.34,-.34);}
\foreach \c in {0,1}{\draw (3.25+.34*\c,.45) rectangle ++(.34,-.34);}
\node at (6.6,.28) {$\varnothing$};
\foreach \x in {9,12}{\draw (\x+.43,.45) rectangle ++(.34,-.34);}
\foreach \x/\sgn in {2.1/{\succ},5.1/{\succ},8.1/{\prec},11.1/{\succ},14.1/{\prec}}{\node at (\x,.15) {$\sgn$};}
\node[anchor=east] at (-.35,-1.05) {$\eta_t$};
\foreach \x/\v in {2.1/0,5.1/0,8.1/1,11.1/0,14.1/1}{\node at (\x,-1.05) {$\v$};}
\end{tikzpicture}
\caption{An example of runner construction for $(a,b)=(2,3)$. Occupied positions are shaded, and the dashed runner is the periodic copy $t=5$ of $t=0$.  The lower rows show staircase removal, conjugation, and the interlacing code. }
\label{fig:residuerunners}
\end{figure}
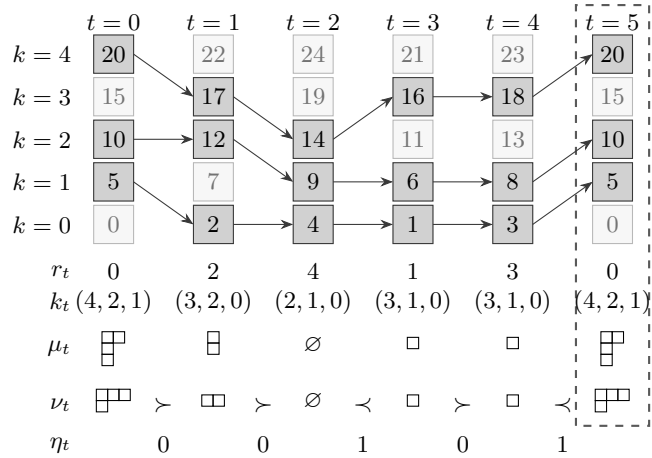

Let the row indices occupied by $\mathcal F$ in the $t$-th column be
\[ k_{t,1}>\dots>k_{t,N}\ge 0.\]
By Lemma~\ref{lem:cyclelength}, the $j$-th highest cycle of $\mathcal F$ is $\dots\to r_t+dk_{t,j}\to r_{t+1}+dk_{t+1,j}\to\dots$.
Treat $k_t=(k_{t,1},\dots,k_{t,N})$ as a partition with distinct parts. Subtracting the staircase gives an ordinary partition with at most $N$ parts:
\[ \mu_t=k_t-(N-1,\dots,0), \quad 0\leq t<d.\]
Equivalently, $\mu_{t,j}$ counts the unoccupied boxes below the $j$-th highest occupied box in the $t$-th runner. Taking conjugates,
\[ \nu_t=\mu_t',\]
gives the sequence in Proposition~\ref{prop:cycleinterlacing}. Its parts are bounded by $N$. Figure~\ref{fig:residuerunners} illustrates the construction.

The construction is reversible. The only remaining requirement is that the edges determined by the occupied vertices be legal. For every $t$ and $j$, this means
\[r_{t+1}+dk_{t+1,j}-r_t-dk_{t,j} \in \{a,-b\}.\]
Since $r_{t+1}-r_t=a-d\eta_t$, this is equivalent to
\begin{equation}\label{eq:runnerdiff}
 k_{t+1,j}-k_{t,j}\in\{\eta_t-1,\eta_t\},
\end{equation}
which is equivalent to
\[ \mu_{t+1,j}-\mu_{t,j}\in \{\eta_t-1,\eta_t\}.\]
If $\eta_t=1$, the condition $\mu_{t+1,j}-\mu_{t,j}\in\{0,1\}$ says that $\mu_{t+1}\supseteq\mu_t$ and that $\mu_{t+1}/\mu_t$ has at most one box in each row. Thus it is a skew column. If $\eta_t=0$, the condition is $\mu_{t+1,j}-\mu_{t,j}\in\{0,-1\}$. Equivalently, $\mu_t\supseteq\mu_{t+1}$ and $\mu_t/\mu_{t+1}$ is a skew column.

Taking the conjugate yields that the corresponding condition on $\nu_t$ is exactly \eqref{eq:interlacing}. This concludes that the proposed process defines a bijection between the sets in Proposition~\ref{prop:cycleinterlacing}. The weight statement \eqref{eq:cyclevolume} follows from \eqref{eq:totalweight} and $\abs{\nu_t}=\abs{\mu_t}=\abs{k_t}-\binom N2$. The proof of Proposition~\ref{prop:cycleinterlacing} is complete.

\subsection{Diagonal slices of cylindric partitions}\label{subsec:diagonalmap}
We now place the interlacing partitions on a cylinder.  The construction is the standard diagonal-slice picture for cylindric partitions \cite[Section~5]{Borodin}.  We record only the coordinates needed to recover our outgoing profile.

Let $s_i,\ i\in \ZZ$ be defined by $s_0=0, s_{i+1}-s_i=c_i$. Draw row $i$ beginning in column $1-s_i$, so consecutive rows are shifted by $c_i$, as is required by the profile. The inequality \eqref{eq:cylinder} is equivalent to the requirement that entries weakly decrease to the east and to the south. For the balanced profile \eqref{eq:balanced}, we have
\[
 s_i:=\left\lfloor\frac{ib}{a}\right\rfloor\qquad(i\in\ZZ).
\]
The picture is periodic under
$(i,\ell)\mapsto(i+a,\ell-b)$, where $\ell$ is the column index.  A northwest--southeast diagonal is called the $t$-th diagonal if its boxes satisfy
\begin{equation}\label{eq:diagonalnumber}
 \ell-i=t+1.
\end{equation}
The first available box on diagonal $t$ occurs in row
\begin{equation}\label{eq:firstrow}
 i_t=-\left\lfloor\frac{at}{d}\right\rfloor.
\end{equation}
Starting there, place the parts of $\nu_t$ successively, one southeast step at a time. The index in $\nu_t$, but not $i_t$, is understood cyclically modulo $d$. The profile determines the starting boxes, while the partitions determine the entries.

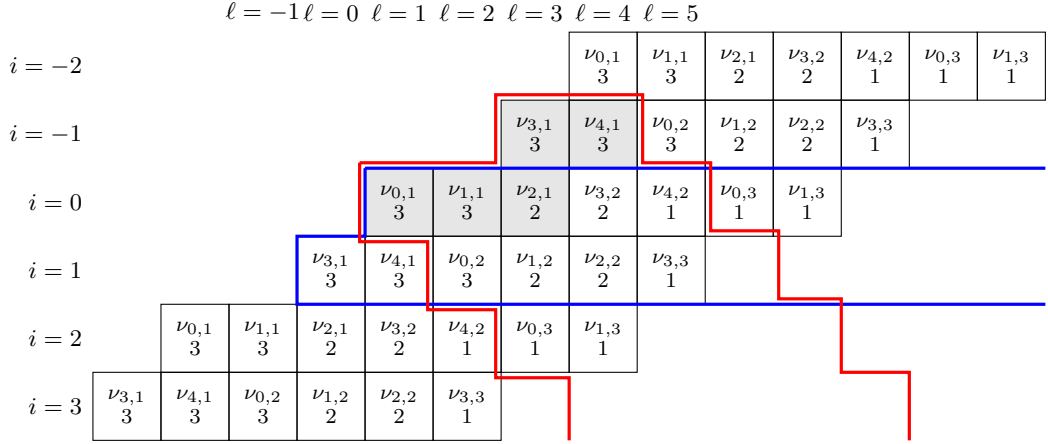
\begin{figure}[!htbp]
\centering
\begin{tikzpicture}[x=.90cm,y=.90cm,font=\tiny]
\node[anchor=east,font=\scriptsize] at (-4,1.50) {$i=-2$};
\node[anchor=east,font=\scriptsize] at (-4,.50) {$i=-1$};
\node[anchor=east,font=\scriptsize] at (-4,-.50) {$i=0$};
\node[anchor=east,font=\scriptsize] at (-4,-1.50) {$i=1$};
\node[anchor=east,font=\scriptsize] at (-4,-2.50) {$i=2$};
\node[anchor=east,font=\scriptsize] at (-4,-3.50) {$i=3$};
\foreach \l in {-1,0,1,2,3,4,5}{
 \node[anchor=south, font=\scriptsize] at (\l-.5,2) {$\ell=\l$};
}
\foreach \i/\l in {0/1,0/2,0/3,-1/3,-1/4}{
\fill[black!10] (\l,-\i) rectangle ++(-1,-1);
}
\foreach \x/\lab/\val in {0/{\nu_{0,1}}/3,1/{\nu_{1,1}}/3,2/{\nu_{2,1}}/2,3/{\nu_{3,2}}/2,4/{\nu_{4,2}}/1,5/{\nu_{0,3}}/1,6/{\nu_{1,3}}/1}{
  \draw (\x,0) rectangle ++(1,-1); \node[align=center] at (\x+.5,-.5) {$\lab$\\[-1pt]$\val$};
\draw (\x+3,+2) rectangle ++(1,-1); \node[align=center] at (\x+.5+3,-.5+2) {$\lab$\\[-1pt]$\val$};
\draw (\x-3,0-2) rectangle ++(1,-1); \node[align=center] at (\x+.5-3,-.5-2) {$\lab$\\[-1pt]$\val$};}
\foreach \x/\lab/\val in {-1/{\nu_{3,1}}/3,0/{\nu_{4,1}}/3,1/{\nu_{0,2}}/3,2/{\nu_{1,2}}/2,3/{\nu_{2,2}}/2,4/{\nu_{3,3}}/1}{
  \draw (\x,-1) rectangle ++(1,-1); \node[align=center] at (\x+.5,-1.5) {$\lab$\\[-1pt]$\val$};
\draw (\x+3,-1+2) rectangle ++(1,-1); \node[align=center] at (\x+.5+3,-1.5+2) {$\lab$\\[-1pt]$\val$};
\draw (\x-3,-1-2) rectangle ++(1,-1); \node[align=center] at (\x+.5-3,-1.5-2) {$\lab$\\[-1pt]$\val$};}
  \draw[blue,very thick] (0,0)--(10,0) (10,-2)--(-1,-2)--(-1,-1)--(0,-1)--(0,0);
\foreach \x in {0.08}{
  \draw[red,very thick] (0-\x,0+\x)--(2-\x,\x)--(2-\x,1+\x)--(4+\x,1+\x)--(4+\x,\x)--(5+\x,\x)--(5+\x,-1+\x)--(6+\x,-1+\x)--(6+\x,-2+\x)--(7,-2+\x)--(7,-3)--(8,-3)--(8,-4) (3,-4)--(3,-3-\x)--(2-\x,-3-\x)--(2-\x,-2-\x)--(1-\x,-2-\x)--(1-\x,-1-\x)--(-\x,-1-\x)--(-\x,\x);}
\end{tikzpicture}
\caption{Diagonal placement for the five slices above. Red region contains the diagonal slices $0$ through $4$. Each $\nu_t$ is placed diagonally, starting from a shaded box. Blue region contains rows $0$ and $1$. Both regions are fundamental domains of the translation $(i,\ell)\mapsto(i+2,\ell-3)$. }
\label{fig:largeslices}
\end{figure}

For example, Figure~\ref{fig:largeslices} shows the placement for $(a,b)=(2,3)$ and
\[
 \nu_0=(3,3,1),\quad \nu_1=(3,2,1),\quad \nu_2=(2,2),\quad
 \nu_3=(3,2,1),\quad \nu_4=(3,1).
\]

The identity $i_{t+1}=i_t-\eta_t$ shows that the interlacing code \eqref{eq:balancedinterlacing} corresponds to the balanced profile in the following sense: adjacent starting boxes differ by exactly the step encoded by $\eta_t$, where $\eta_t=0$ corresponds to an east step from $t$-th starting box to the $(t+1)$-th, and $\eta_t=1$ a north step.  The interlacing conditions in \eqref{eq:interlacing}, interpreted by the inequalities \eqref{eq:interlacing-inequality}, are therefore equivalent to weak decrease along rows and columns on the cylinder. This proves that the placement is bijective.

The periodic translation identifies diagonal $t+d$ with diagonal $t$.  Thus the diagonals $0,\ldots,d-1$ contain one representative of every box in a fundamental cylinder domain, as do rows $0,\dots,a-1$. Since the placement does not change any entry,
\begin{equation}\label{eq:volume}
 |\boldsymbol\lambda|=\sum_{t=0}^{d-1}|\nu_t|.
\end{equation}
The largest cylinder entry is the largest part among the slices, so it is at most $N$.  Together with Proposition~\ref{prop:cycleinterlacing} and \eqref{eq:cyclevolume}, this proves Proposition~\ref{prop:cyclecylinder}.

\begin{remark}[The zero cylinder]
For any $N$, the zero cylinder in Proposition~\ref{prop:cyclecylinder} gives empty slices $\nu_t$ and $\mu_t$, hence the staircase levels $k_t=(N-1,\dots,0)$ in every runner. The corresponding cycle collection occupies every vertex in $[0,dN-1]$. By \eqref{eq:cyclevolume}, its weight is $\gamma_N$. Thus the zero cylinder already accounts for the shift in \eqref{eq:cyclecylindertarget}.
\end{remark}

\section{The finite identity and its consequences}\label{sec:conclusion}
We combine the two solutions of the common $q$-difference equation
to prove the finite identity.  We then record the four equivalent
generating functions and take the limit that proves the HJO conjecture.

\begin{proof}[Proof of Theorems~\ref{thm:common} and~\ref{thm:mainfinite}]
Equation~\eqref{eq:HP} gives $\Hseries=\Pseries$, and
Proposition~\ref{prop:cyclecylinder} gives $\Bseries=\Dseries$.
By Theorems~\ref{thm:Euler} and~\ref{thm:cofactor}, both series
satisfy \eqref{eq:commonequation}.  They have constant coefficient
$1$ and belong to $\QQ((q))[[z]]$.  The uniqueness established by
\eqref{eq:solutionrecurrence} therefore gives
\[
 \Hseries=\Pseries=\Dseries=\Bseries
      =\sum_{N\ge0}G_N(q)z^N.
\]
Their coefficients satisfy
\[
 q^{\gamma_N}\frac{F_N(q)}{(q)_N}
       =q^{\gamma_N}C_{\cc,\le N}(q)=G_N(q).
\]
Multiplying by $q^{-\gamma_N}(q)_N$ proves
\eqref{eq:commoncoefficients} and \eqref{eq:mainfinite}.
\end{proof}

For every $N\ge 0$, recall the balanced profile \eqref{eq:balanced}, the generating function $P_N(q)$ from \eqref{eq:pathseries}, and $D_N(q)$ from \eqref{eq:cycleseries}. The finite conclusion of this paper can be succinctly summarized as follows: three normalized combinatorial generating functions are polynomials with constant term $1$, and they are equal. Moreover, they are given by a common algebraic formula $F_N(q)$.

\begin{corollary}\label{cor:threemodels}
For every $N\ge 0$ and coprime $a,b>1$, with the balanced profile $\cc$ from \eqref{eq:balanced}, we have
\begin{equation}\label{eq:threeway}
 \begin{aligned}
 F_N(q)&=(q;q)_N C_{\cc,\le N}(q)\\
       &=q^{-\kappa_N}P_N(q)
        =q^{-\gamma_N}(q;q)_N D_N(q).
 \end{aligned}
\end{equation}
\end{corollary}
\begin{proof}
Combine Theorem~\ref{thm:mainfinite}, Corollary~\ref{prop:return},
and Proposition~\ref{prop:cyclecylinder}.
\end{proof}

The cycle and cylinder descriptions are related by a bijection.
Their equality with the return-path series follows from the common
equation, not from a direct bijection between paths and cylinders.

Finally, the infinite identity conjectured by HJO follows by taking coefficientwise limits and
using the known cylinder product.

\begin{proof}[Proof of Corollary~\ref{cor:HJO}]
For any fixed power $q^j$, bound \eqref{eq:coercivity} leaves only
finitely many gap vectors that can contribute.  For each such vector,
$(q)_N/(q)_{N-n_f}$ tends coefficientwise to $1$, so
$F_N(q)\to Z_{a,b}(q)$.  A cylinder of volume at most $j$ has every
entry at most $j$, hence $C_{\cc,\le N}(q)\to C_\cc(q)$.
Taking limits in \eqref{eq:mainfinite} gives
\[
 Z_{a,b}(q)=(q;q)_\infty C_\cc(q)=P_{a,b}(q).
\]
The last equality is the balanced-profile product of
\cite[Proposition~9.3]{HJO} and \cite[Sections~2--3]{FW}, as recalled
in the introduction.
\end{proof}

\section{\texorpdfstring{A further finite example}{A further finite example}}\label{sec:examples}
We illustrate the finite identity for $(a,b)=(4,5)$ by computing the
first two HJO polynomials.  These determine the three combinatorial generating functions in Corollary~\ref{cor:threemodels}.

For $(a,b)=(4,5)$, the gap set and balanced profile are
\[
 G=\{1,2,3,6,7,11\},\qquad \delta=6,\qquad \cc=(1,1,1,2).
\]
The number of order filters is $14$, the fourth Catalan number. Their rank-one weights give
\[
 F_1(q)=1+3q+3q^2+3q^3+2q^4+q^5+q^6.
\]
At rank two the original Gaussian sum yields
\begin{align*}
 F_2(q)={}&1+3q+6q^2+9q^3+13q^4+14q^5+17q^6+16q^7\\
 &+17q^8+15q^9+15q^{10}+12q^{11}+13q^{12}+9q^{13}+9q^{14}\\
 &+6q^{15}+6q^{16}+4q^{17}+4q^{18}+2q^{19}+2q^{20}\\
 &+q^{21}+q^{22}+q^{24}.
\end{align*}
We have $F_2(1)=196=14^2$, as the filter-word formula predicts, and Corollary~\ref{cor:threemodels} gives
\[
 C_{(1,1,1,2),\le2}(q)=\frac{F_2(q)}{(1-q)(1-q^2)},
\]
\[
D_2(q)=q^9\frac{F_2(q)}{(1-q)(1-q^2)},
\]
and
\[ P_2(q)=q^8F_2(q).\]
The two shifts are $\kappa_2=8$ and $\gamma_2=9$.

\section{AxiomProver Lean certificate}\label{sec:AxiomProver}

This paper was written using human-AI collaboration with AxiomProver. This paper is the third in a series, and it proves the full HJO conjecture. The previous two articles \cite{HLOP, LO} offered some evidence for the conjecture, with proofs involving more standard techniques. For the present paper, the human-AI collaboration built on these earlier two papers, and proceeded through literature search, significant numerical experimentation, and the proposal and rejection of many candidate strategies before the proofs below took shape. At every stage the authors set the direction and made the mathematical judgments, and they take full responsibility for the correctness of the results. AxiomProver is an AI system for mathematical research through formal proof, currently under development at Axiom Math.

 The finite identity of
Theorem~\ref{thm:mainfinite} and the HJO conjecture in
Corollary~\ref{cor:HJO} have been formalized with AxiomProver in Lean. The repository contains
the formal statements, proofs, and instructions for checking the certificate:
\begin{center}
\url{https://github.com/AxiomMath/HJO}
\end{center}
The human-written challenge file gives the definitions and hypotheses.
The solution proves the stated targets. The repository reports
successful verification with the Lean comparator.

\end{document}